\documentclass{amsproc}

\usepackage{amsmath,amssymb}
\newtheorem{theorem}{Theorem}[section]
\newtheorem{lemm}[theorem]{Lemma}
\newtheorem{prop}[theorem]{Proposition}

\theoremstyle{definition}
\newtheorem{defi}[theorem]{Definition}

\theoremstyle{remark}
\newtheorem{remark}[theorem]{Remark}

\numberwithin{equation}{section}

\def\lg{\langle}
\def\rg{\rangle}

\def\vn{\varepsilon}

\def\ot{\otimes}

\def\om{\omega}

\newfont{\df}{eufm10}

\def\ot{\otimes}

\def\ot{\otimes}

\begin{document}

\title[Two-parameter quantum group of exceptional type $G_2$]
{Two-parameter Quantum Group of Exceptional Type $G_2$\\ and
Lusztig's Symmetries}

\author[Hu]{Naihong Hu$^\star$}
\address{Department of Mathematics, East China Normal University,
Shanghai 200062, PR China} \email{nhhu@euler.math.ecnu.edu.cn}
\thanks{$^\star$N.H., Corresponding Author,
supported in part by the NNSF (Grant 10431040), the TRAPOYT,
the FUDP and the PCSIRT from the MOE of China, the SRSTP from the STCSM
}

\author[Shi]{Qian Shi}
\address{Department of Mathematics, East China Normal University,
Shanghai 200062, PR China}

\date{June 18, 2005}


\begin{abstract}
We give the defining structure of two-parameter quantum group of
type $G_2$ defined over a field ${\Bbb Q}(r,s)$ (with $r\ne s$), and
prove the Drinfel'd double structure as its upper and lower
triangular parts, extending an earlier result of [BW1] in type $A$
and a recent result of [BGH1] in types $B, C, D$. We further discuss
the Lusztig's $\Bbb Q$-isomorphisms from $U_{r,s}(G_2)$ to its
associated object $U_{s^{-1},r^{-1}}(G_2)$, which give rise to the
usual Lusztig's symmetries defined not only on $U_q(G_2)$ but also
on the centralized quantum group $U_q^c(G_2)$ only when
$r=s^{-1}=q$. (This also reflects the distinguishing difference
between our newly defined two-parameter object and the standard
Drinfel'd-Jimbo quantum groups). Some interesting $(r,s)$-identities
holding in $U_{r,s}(G_2)$ are derived from this discussion.
\end{abstract}

\maketitle

\section{Two-parameter Quantum Group $U_{r,s}(G_2)$} 
\medskip

Let ${\Bbb K}={\Bbb Q}(r,s)$ be a field of rational functions with
two indeterminates $r$, $s$.

Let $\Phi$ be a finite root system of $G_2$ with $\Pi$ a base of
simple roots, which is a subset of a Euclidean space $E = {\Bbb
R}^3$ with an inner product $(\,,\,)$. Let $\epsilon
_{1},\,\epsilon _{2},\,\epsilon_{3}$ denote an orthonormal basis
of $E$, then $\Pi = \{\alpha_{1} = \epsilon_{1}-\epsilon_{2},\;
\alpha_{2} = \epsilon_{2}+\epsilon_{3}-2\epsilon_{1}\}$ and $\Phi
= \pm
\{\alpha_1,\alpha_2,\alpha_2+\alpha_1,\alpha_2+2\alpha_1,\alpha_2+3\alpha_1,2\alpha_2+3\alpha_1\}$.
In this case, we set $\displaystyle r_1 =
r^{\frac{(\alpha_1,\,\alpha_1)}{2}} = r,\; r_2 =
r^{\frac{(\alpha_2,\,\alpha_2)}{2}} = r^3$ and $s_1 =
s^{\frac{(\alpha_1,\,\alpha_1)}{2}} = s,\; s_2 =
s^{\frac{(\alpha_2,\,\alpha_2)}{2}} = s^3$.

We begin by giving the definition of two-parameter quantum group
of type $G_2$, which is new.

\begin{defi}
Let $U=U_{r,\,s}(G_2)$ be the associative algebra over ${\Bbb
Q}(r,s)$ generated by symbols $e_i,\;f_i,\;\omega_i^{\pm 1},\;
\omega_i'^{\pm 1} \;(1\leq i\leq 2)$ subject to the relations
\vskip0.2cm

\noindent $(G1)$ \ $[\,\omega_i^{\pm 1}, \omega_j^{\pm
1}\,]=[\,\omega_i^{\pm 1}, \omega_j'^{\pm 1}\,]=[\,\omega_i'^{\pm
1}, \omega_j'^{\pm 1}\,]=0, \quad \omega_i\omega_i^{-1} =1=
\omega_j'\omega_j'^{-1}$.

\noindent $(G2)$
\hspace*{\fill}$\omega_{1}\,e_{1}\,\omega_{1}^{-1} =
(rs^{-1})\,e_{1},\qquad\quad\ \;
\omega_{1}\,f_{1}\,\omega_{1}^{-1} =
(r^{-1}s)\,f_{1}$,\hspace*{\fill}

 \hspace*{\fill}$\omega_{1}\,e_{2}\,\omega_{1}^{-1} =
s^{3}\,e_{2},\qquad\qquad\quad \omega_{1}\,f_{2}\,\omega_{1}^{-1}
= s^{-3}\,f_{2}$,\hspace*{\fill}

\hspace*{\fill}$\omega_{2}\,e_{1}\,\omega_{2}^{-1} =
r^{-3}\,e_{1},\qquad\quad\ \;\, \omega_{2}\,f_{1}\,\omega_{2}^{-1}
= r^{3}\,f_{1}$,\hspace*{\fill}

\hspace*{\fill}$\omega_{2}\,e_{2}\,\omega_{2}^{-1} =
(r^{3}s^{-3})\,e_{2},\qquad\qquad\
\omega_{2}\,f_{2}\,\omega_{2}^{-1} =
(r^{-3}s^{3})\,f_{2}$.\hspace*{\fill}

\

\noindent $(G3)$
\hspace*{\fill}$\omega_{1}'\,e_{1}\,\omega_{1}'^{-1} =
(r^{-1}s)\,e_{1},\qquad\quad\ \
\omega_{1}'\,f_{1}\,\omega_{1}'^{-1} =
(rs^{-1})\,f_{1}$,\hspace*{\fill}

 \hspace*{\fill}$\omega_{1}'\,e_{2}\,\omega_{1}'^{-1} =
r^{3}\,e_{2},\qquad\qquad\quad
\omega_{1}'\,f_{2}\,\omega_{1}'^{-1} =
r^{-3}\,f_{2}$,\hspace*{\fill}

\hspace*{\fill}$\omega_{2}'\,e_{1}\,\omega_{2}'^{-1} =
s^{-3}\,e_{1},\qquad\quad\ \;\,
\omega_{2}'\,f_{1}\,\omega_{2}'^{-1} =
s^{3}\,f_{1}$,\hspace*{\fill}

\hspace*{\fill}$\omega_{2}'\,e_{2}\,\omega_{2}'^{-1} =
(r^{-3}s^{3})\,e_{2},\qquad\qquad\
\omega_{2}'\,f_{2}\,\omega_{2}'^{-1} =
(r^{3}s^{-3})\,f_{2}$.\hspace*{\fill}

\

\noindent $(G4)$ \ For $1\le i,\, j\le 2$, we have
$$
[\,e_i, f_j\,]=\delta_{ij}\frac{\om_i-\om_i'}{r_i-s_i}.
$$

\noindent $(G5)$ \, ($(r,\,s)$-Serre relations)
\begin{gather*}
e_2^{2}e_1 - (r^{-3} + s^{-3})\,e_{2}e_{1}e_{2} +
(rs)^{-3}\,e_{1}e_2^2 = 0,\tag*{$(G5)_1$}\\
\begin{split}
e_1^{4}e_2 - (r +& s)(r^2 + s^2)\,e_1^{3}e_{2}e_1 + rs(r^2 +
s^2)(r^2
+ rs + s^2)\,e_{1}^{2}e_{2}e_1^2\\
-& (rs)^3(r + s)(r^{2} + s^2)\,e_{1}e_{2}e_1^3+\, (rs)^6e_{2}e_1^4
= 0.
\end{split}\tag*{$(G5)_2$}
\end{gather*}

\noindent $(G6)$ \, ($(r,\,s)$-Serre relations)
\begin{gather*}
f_{1}f_2^{2} - (r^{-3} + s^{-3})\,f_{2}f_{1}f_{2} +
(rs)^{-3}\,f_{2}^{2}f_1 = 0,\tag*{$(G6)_1$}\\
\begin{split}
f_{2}f_1^{4} - (r +& s)(r^2 + s^2)\,f_1f_{2}f_1^{3} + rs(r^2 +
s^2)(r^2 + rs + s^2)\,f_{1}^{2}f_{2}f_1^2\\
-& (rs)^3(r + s)(r^{2} + s^2)\,f_{1}^3f_{2}f_1 +\,
(rs)^6\,f_{1}^4f_2 = 0.
\end{split}\tag*{$(G6)_2$}
\end{gather*}

\end{defi}

\begin{prop}
The algebra $U_{r,s}(G_2)$ is a Hopf algebra with comultiplication,
counit and antipode given by
\begin{gather*}
\Delta(\om_i^{\pm1})=\om_i^{\pm1}\ot\om_i^{\pm1}, \qquad
\Delta({\om_i'}^{\pm1})={\om_i'}^{\pm1}\ot{\om_i'}^{\pm1},\\
\Delta(e_i)=e_i\ot 1+\om_i\ot e_i, \qquad \Delta(f_i)=1\ot
f_i+f_i\ot \om_i',\\
\vn(\om_i^{\pm})=\vn({\om_i'}^{\pm1})=1, \qquad
\vn(e_i)=\vn(f_i)=0,\\
S(\om_i^{\pm1})=\om_i^{\mp1}, \qquad
S({\om_i'}^{\pm1})={\om_i'}^{\mp1},\\
S(e_i)=-\om_i^{-1}e_i,\qquad S(f_i)=-f_i\,{\om_i'}^{-1}.
\end{gather*}\hfill\qed
\end{prop}

\begin{remark} \ (I) \ When $r = q=s^{-1}$, the quotient Hopf algebra of
$U_{r,\,s}(G_2)$ modulo the Hopf ideal generated by elements
$\omega_i' - \omega_i^{-1} (1\leq i \leq 2)$, is just the standard
quantum group $U_q(G_2)$ of Drinfel'd-Jimbo type; the quotient
modulo the Hopf ideal generated by elements $\om_i'-z_i\om_i^{-1}$
$(1\le i\le 2)$, where $z_i$ runs over the center, is the {\it
centralized quantum group} $U_q^c(G_2)$.

(II) \ In any Hopf algebra $\mathcal H$, there exist the
left-adjoint and the right-adjoint action defined by the Hopf
algebra structure:
$$
\text{ad}_{ l}\,a\,(b)=\sum_{(a)}a_{(1)}\,b\,S(a_{(2)}), \qquad
\text{ad}_{ r}\,a\,(b)=\sum_{(a)}S(a_{(1)})\,b\,a_{(2)},
$$
where $\Delta(a)=\sum_{(a)}a_{(1)}\ot a_{(2)}\in \mathcal H\otimes
\mathcal H$, for any $a$, $b\in \mathcal H$.

From the viewpoint of adjoint actions, the $(r,s)$-Serre relations
$(G5)$, $(G6)$ take on the simpler forms
\begin{gather*}
\bigl(\text{ad}_l\,e_i\bigr)^{1-a_{ij}}\,(e_j)=0,
\qquad\text{\it for any } \ i\ne j,\\
\bigl(\text{ad}_r\,f_i\bigr)^{1-a_{ij}}\,(f_j)=0, \qquad\text{\it
for any } \ i\ne j.
\end{gather*}

\end{remark}

\bigskip
\section{Drinfel'd Quantum Double}
\medskip

\begin{defi}
A (Hopf) dual pairing of two Hopf algebras $\mathcal{A}$ and
$\mathcal{U}$ (see [BGH1], or [KS]) is a bilinear form $\langle \,
, \rangle:\, \mathcal{U} \times \mathcal{A} \longrightarrow
\mathbb{K}$ such that
\begin{gather*}
\langle f,\,1_\mathcal{A} \rangle =
\varepsilon_\mathcal{U}(f),\qquad \langle 1_\mathcal{U},\, a
\rangle = \varepsilon_\mathcal{A}(a),\tag{1}\\
\langle f,\, a_{1}a_2 \rangle = \langle
\triangle_\mathcal{U}(f),\,a_1\otimes a_2 \rangle, \qquad \langle
f_1{f}_2,\,a \rangle = \langle f_1\otimes
f_2,\,\triangle_\mathcal{A}(a) \rangle,\tag{2}
\end{gather*}
for all $f,\, f_1,\, f_2 \in \mathcal{U}$, and $a,\, a_1,\, a_2
\in \mathcal{A}$, where $\varepsilon_\mathcal{U}$ and
$\varepsilon_\mathcal{A}$ denote the counits  of $\mathcal{U}$ and
$\mathcal{A}$, respectively, and $\triangle_\mathcal{U}$ and
$\triangle_\mathcal{A}$ are their comultiplications.
\end{defi}

A direct consequence of the defining properties above is that
$$\langle S_\mathcal{U}(f),\,a \rangle = \langle
f,\,S_\mathcal{A}(a) \rangle, \quad f\in \mathcal{U},\,a\in
\mathcal{A},$$
where $S_\mathcal{U}, S_\mathcal{A}$ denote the
respective antipodes of $\mathcal{U}$ and $\mathcal{A}$.

\begin{defi}
A bilinear form $\langle \,,\, \rangle:\,\mathcal{U} \times
\mathcal{A} \longrightarrow \mathbb{K}$ is called a skew-dual
pairing of two Hopf algebras $\mathcal{A}$ and $\mathcal{U}$ (see
[BGH1]) if $\langle\, ,\rangle: \, \mathcal{U}^{\mathrm {cop}}
\times \mathcal{A} \longrightarrow \mathbb{K}$ is a Hopf dual
pairing of $\mathcal{A}$ and $\mathcal{U}^{\mathrm {cop}}$, where
$\mathcal{U}^{\mathrm {cop}}$ is the Hopf algebra having the
opposite comultiplication to $\mathcal{U}$, and
$S_{\mathcal{U}^{\mathrm {cop}}} = S_\mathcal{U}^{-1}$ is
invertible.
\end{defi}

Denote by $\mathcal{B} = B(G_2)$   the Hopf subalgebra of
$U_{r,\,s}(G_2)$ generated by $e_j,\omega_j^{\pm 1}$ and by
$\mathcal{B'} = B'(G_2)$ the one generated by $f_j, \omega_j'^{\pm
1}$, where $1 \leq j\leq 2$.

\begin{prop} There exists a unique
skew-dual pairing $\langle \,, \rangle: \,\mathcal{B'} \times
\mathcal{B} \longrightarrow \Bbb Q(r,\,s)$ of the Hopf subalgebras
$\mathcal{B}$ and $\mathcal{B'}$, such that
\begin{gather*}
\langle f_i,\,e_j \rangle = \delta_{ij} \frac{1}{s_i - r_i},
\qquad (1 \leq i,\,j \leq2),\tag{3}\\
\langle \omega_1',\,\omega_1 \rangle = rs^{-1},\tag{$4_1$}\\
\langle \omega_1',\,\omega_2 \rangle = r^{-3}, \tag{$4_2$} \\
\langle \omega_2',\,\omega_1 \rangle = s^{3}, \tag{$4_3$}\\
\langle \omega_2',\,\omega_2 \rangle = r^{3}s^{-3},\tag{$4_4$}\\
\langle \omega_i'^{\pm 1},\,\omega_j^{-1} \rangle = \langle
\omega_i'^{\pm 1},\,\omega_j \rangle ^{-1} = \langle
\omega_i',\,\omega_j \rangle ^{\mp 1},\qquad (1 \leq i,\,j \leq
2),\tag{5}
\end{gather*}
and all other pairs of generators are $0$. Moreover, we have
$\langle S(a),\,S(b) \rangle = \langle a,\,b \rangle$ for $a\in
\mathcal{B}', b\in \mathcal{B}$.
\end{prop}
\begin{proof}  \  Since any skew-dual pairing of bialgebras is
determined by its values on generators, uniqueness is clear. We
proceed to prove the existence of the pairing.

We begin by defining a bilinear form $\langle\,, \rangle :
\mathcal{B}'^{\mathrm {cop}} \times \mathcal{B} \rightarrow {\Bbb
Q}(r,\,s)$ first on the generators satisfying (3), (4), and (5).
Then we extend it to a bilinear form on $ \mathcal{B}'^{\mathrm
{cop}} \times \mathcal{B}$ by requiring that (1) and (2) hold for
$\triangle_{\mathcal{B}'^{\mathrm {cop}}} =
\triangle_{\mathcal{B}'}^{\mathrm {op}}$. We will verify that the
relations in $\mathcal{B}$ and $\mathcal{B}'$ are preserved,
ensuring that the form is well-defined and so is a dual pairing of
$\mathcal{B}$ and $\mathcal{B}'^{\mathrm {cop}}$ by definition.

It is direct to check that the bilinear form preserves all the
relations among the $\omega_i^{\pm 1}$ in $\mathcal{B}$ and the
$\omega_i'^{\pm 1}$ in $\mathcal{B}'$. Next, the structure constants
$(4_n)$ ensure the compatibility of the form defined above with
those relations of $(G2)$ and $(G3)$ in $\mathcal{B}$ or
$\mathcal{B}'$ respectively. We are left to verify that the form
preserves the $(r,\,s)$-Serre relations in $\mathcal{B}$ and
$\mathcal{B}'$. It suffices to show that the form on
$\mathcal{B}'^{\mathrm {cop}} \times \mathcal{B}$ preserves the
$(r,s)$-Serre relations in $\mathcal{B}$; the verification for
$\mathcal{B}'^{\mathrm {cop}}$ is similar.

\smallskip
First, let us show that the form preserves the $(r,\,s)$-Serre
relation of degree $2$ in $\mathcal{B}$, that is,
$$
\langle X,\,e_2^{2}e_1 - (r^{-3} + s^{-3})\,e_{2}e_{1}e_{2} +
r^{-3}s^{-3}\,e_{1}e_2^2 \,\rangle = 0,
$$
where $X$ is any word in the generators of $\mathcal{B}'$. It
suffices to consider three monomials: $X = f_2^{2}f_1,\,
f_{2}f_{1}f_2,\, f_{1}f_2^2$. However, in the degree $2$'s
situation for type $G_2$, its proof is the same as that of type
$C_2$ (see [BGH1, (7C) and thereafter]).

\smallskip
Next, we verify that the $(r,\,s)$-Serre relation of degree $4$ in
$\mathcal{B}$ is preserved by the form, that is, we show that
\begin{gather*}
\langle X, \; e_1^{4}e_2 - (r + s)(r^2 + s^2)\,e_1^{3}e_{2}e_1 +
rs(r^2 + s^2)(r^2 + rs + s^2)\,e_{1}^{2}e_{2}e_1^2\\
- (rs)^3(r + s)(r^{2} + s^2)\,e_{1}e_{2}e_1^3+\, (rs)^6\,e_{2}e_1^4
\, \rangle,
\end{gather*}
vanishes, where $X$ is any word in the generators of $\mathcal{B}'$.
By definition, this expression equals
\begin{equation*}
\begin{split}
\langle \triangle^{(4)}(X),\;&\,e_1 \otimes e_1 \otimes e_1
\otimes
e_1 \otimes e_2 \\
&- (r + s)(r^2 + s^2)\,e_1 \otimes e_1 \otimes e_1
\otimes e_{2}\otimes e_1\\
& + rs(r^2 + s^2)(r^2 + rs + s^2)\,e_{1}\otimes e_{1} \otimes e_{2
} \otimes e_1 \otimes e_{1} \\
& - (rs)^3(r + s)(r^{2} + s^2)\,e_{1} \otimes e_{2} \otimes e_1
\otimes e_1 \otimes e_1\\
& + (rs)^6\,e_{2} \otimes e_1 \otimes e_1 \otimes e_1 \otimes e_1 \,
\rangle,
\end{split}\tag{$\ast$}
\end{equation*}
where $\triangle$ in the left-hand side of the pairing
$\langle\,,\rangle$ indicates $\triangle_{\mathcal{B}'}^{\mathrm
{op}}$. In order for any one of these terms to be nonzero, $X$
must involve exactly four $f_1$ factors, one $f_2$ factor, and
arbitrarily many $\omega_j'^{\pm 1}$ factors $(j = 1, 2)$.

\smallskip
It suffices to consider five key cases:

(i) $X = f_1^{4}f_2$, we have
\begin{equation*}
\begin{split}
\triangle^{(4)}(X)& = \bigl(\,\omega_1' \otimes \omega_1' \otimes
\omega_1' \otimes \omega_1' \otimes f_1 + \omega_1' \otimes
\omega_1' \otimes \omega_1' \otimes f_1 \otimes 1\\
&\quad+ \omega_1' \otimes \omega_1' \otimes f_1 \otimes 1 \otimes
1 +\omega_1'\otimes f_1 \otimes 1 \otimes 1 \otimes 1 + f_1
\otimes
1 \otimes 1 \otimes 1 \otimes 1 \,\bigr)^4 \cdot \\
&\quad\bigl(\,\omega_2' \otimes \omega_2' \otimes \omega_2'
\otimes \omega_2' \otimes f_2 + \omega_2' \otimes \omega_2'
\otimes \omega_2' \otimes f_2 \otimes
1\\
&\quad+ \omega_2' \otimes \omega_2' \otimes f_2 \otimes 1 \otimes
1 + \omega_2' \otimes f_2 \otimes 1 \otimes 1 \otimes 1 + f_2
\otimes 1 \otimes 1 \otimes 1 \otimes 1\,\bigr).
\end{split}
\end{equation*}
Expanding $\triangle ^{(4)}(X)$, we get $120$ relevant terms having
a nonzero contribution to $(\ast)$. They are listed in TABULAR $1$
of Appendix, together with their pairing values, where we have
introduced
$$a = \langle f_1,\,e_1 \rangle^4 \langle f_2,\, e_2
\rangle, \quad x = \langle \omega_1',\,\omega_1 \rangle, \quad \bar
x = \langle \omega_1',\,\omega_2 \rangle,\quad y=\langle
\omega_2',\,\omega_1 \rangle.$$

The expression in $(\ast)$ equals
\begin{equation*}
\begin{split}
&(\hbox{sum of expressions in column 1})\\
& - (\text{sum of expressions in column 2}) \cdot (r + s)(r^2
+ s^2)\\
& + (\hbox{sum of expressions in column 3}) \cdot rs(r^2 +
s^2)(r^2 + rs + s^2)\\
& - (\hbox{sum of expressions in column 4}) \cdot (rs)^3(r +
s)(r^2 + s^2)\\
& + (\hbox{sum of expressions in column 5}) \cdot (rs)^6.
\end{split}
\end{equation*}

Thus, if we sum up all the pairing values listed in each column of
TABULAR $1$ and multiply by the appropriate factor, we obtain the
paring value of ($\ast$):
\begin{equation*}
\begin{split}
& a(1 + 3x + 5x^2 + 6x^3 + 5x^4 + 3x^5 + x^6)\cdot[\,1 - (r +
s)(r^2
+ s^2)\bar x\\
&\quad + rs(r^2 + s^2)\cdot (r^2 + rs +s^2)\bar x^2 - (rs)^3(r +
s)(r^2 + s^2)\bar x^3 +
(rs)^6 \bar x^4\,]\\
&= a(1 + 3x + 5x^2 + 6x^3 + 5x^4 + 3x^5 + x^6)(1{-}r^{3}\bar x)(1
{-}r^{2}s\bar x)(1{-}rs^{2}\bar x)(1{-}s^{3}\bar x )\\
&= 0 \qquad\qquad \hbox{(because $\bar x = \langle
\omega_1',\,\omega_2 \rangle = r^{-3}$\,).}
\end{split}
\end{equation*}

(ii) $X = f_{2}f_1^4$. By calculation, we get $120$ relevant terms
of $\triangle ^{(4)}(X)$ in ($\ast$) and their pairing values listed
in TABULAR $2$ of Appendix.

If we sum up all the pairing values listed in each column of TABULAR
$2$, then we obtain the paring values of ($\ast$):
\begin{equation*}
\begin{split}
&a(1 + 3x + 5x^2 + 6x^3 + 5x^4 + 3x^5 + x^6)\cdot[\,y^4 - (s^{3} +
rs^{2} + r^{2}s + r^3)\cdot y^3\\
&\quad + rs(s^4 + rs^3 + 2r^{2}s^2 + r^{3}s + r^4)\cdot y^2 -
(rs)^3(s^{3} + rs^{2} + r^{2}s +
r^3)\cdot y + (rs)^6\,]\\
&= a(1 + 3x + 5x^2 + 6x^3 + 5x^4 + 3x^5 + x^6)(y -
r^3)(y - r^{2}s)(y - rs^2)(y - s^3)\\
&= 0 \qquad\qquad (\text{because } y = \langle \omega_2',\,\omega_1
\rangle = s^{3}).
\end{split}
\end{equation*}

(iii) $X = f_{1}^{2}f_{2}f_{1}^{2}$. By calculation, we get $120$
relevant of $\triangle ^{(4)}(X)$ in ($\ast$) and their pairing
values listed in TABULAR $3$ of Appendix.

If we sum up all the pairing values listed in each column of TABULAR
$3$, then we get the paring values of ($\ast$):
\begin{equation*}
\begin{split}
& ay^{2}(1 + 3x + 5x^2 + 6x^3 + 5x^4 + 3x^5 + x^6) \\
&\quad- ay(s^3 + rs^2 + r^{2}s + r^3)\cdot(1 + 3x + 4x^2 + 3x^3 +
x^4)\cdot(1 +
\bar{x}yx^2)\\
&\quad + ars(s^4 + rs^3 + 2r^{2}s^2 +r^{3}s + r^4)\cdot[\,1 + 2x + x^2\\
&\quad + \bar{x}xy(1 + 4x + 6x^2 +4x^3 + x^4) +
\bar{x}^{2}y^{2}x^4\cdot(1+ 2x + x^2)\,]\\
&\quad - a(rs)^3(s^3 + rs^2 + r^{2}s + r^3)\cdot(1 + 3x + 4x^2 +
3x^3
+ x^4)\cdot(\bar{x} + \bar{x}^{2}x^2y)\\
&\quad + (rs)^6a\bar{x}^2( 1 + 3x + 5x^2 + 6x^3 + 5x^4 + 3x^5 + x^6)
\\
&= 2ar^{-1}s^{-1}\cdot rs\cdot (s^4 + 3rs^3 + 4r^{2}s^2 + 3r^{3}s
+ r^4)\cdot(s^2 + r^2)\\
&\quad +  2ar^{-1}s^{-1}\cdot (s^4 + 3rs^3 + 4r^{2}s^2 + 3r^{3}s +
r^4)\cdot (s^2 + rs + r^{2})\cdot(s^2 + r^2)\\
&\quad - 2ar^{-1}s^{-1}\cdot(s^4 + 3rs^3 + 4r^{2}s^2 + 3r^{3}s +
r^4)\cdot(s^2 + r^2)\cdot(r + s)^2\\
&= 2ar^{-1}s^{-1}\cdot (s^4 + 3rs^3 + 4r^{2}s^2 + 3r^{3}s +
r^4)\cdot(s^2 + r^2)\cdot(r + s)^2 \\
&\quad - 2ar^{-1}s^{-1}\cdot(s^4 +
3rs^3 + 4r^{2}s^2 + 3r^{3}s + r^4)\cdot(s^2 + r^2)\cdot(r + s)^2\\
&= 0 \qquad (\text{because } x = \langle \omega_1',\,\omega_1
\rangle = rs^{-1},\, \bar{x} = \langle \omega_1',\,\omega_2 \rangle
= r^{-3},\, y = \langle \omega_2',\,\omega_1 \rangle = s^{3}).
\end{split}
\end{equation*}

(iv) $X = f_{1}^{3}f_{2}f_{1}$. By calculation, we get $120$
relevant terms of $\triangle ^{(4)}(X)$ in ($\ast$) and their
pairing values listed in TABULAR $4$ of Appendix.

If we sum up all the pairing values listed in each column of TABULAR
$4$, then we get the paring-values of ($\ast$):
\begin{equation*}
\begin{split}
& ay(1 {+} 3x {+} 5x^2 {+} 6x^3 {+} 5x^4 {+} 3x^5 {+} x^6) - a(r {+}
s)(r^2 {+} s^2)[\,1 {+} 2x {+} 2x^2 {+} x^3\\
&\quad + x\bar{x}y(1 {+} 3x {+} 5x^2 {+} 5x^3 {+} 3x^4 {+} x^5)\,]+
ars(r^2 {+}
s^2)(r^2 {+} rs {+} s^2)\cdot\\
&\quad[\,\bar{x}(1 {+} 3x {+} 4x^2 {+}
3x^3 {+} x^4)+ \bar{x}^{2}yx^2(1 {+} 3x {+} 4x^2 {+} 3x^3 {+} x^4)\,]\\
&\quad - a(rs)^3(r {+} s)(r^2 {+} s^2)[\,\bar{x}^{2}(1 {+} 3x {+}
5x^2 {+}
5x^3 {+} 3x^4 {+} x^5)\\
&\quad + \bar{x}^{3}y(x^3 {+} 2x^4 {+} 2x^5 {+} x^6)\,]+
ar^{6}s^6\bar{x}^3(1 {+} 3x {+} 5x^2 {+} 6x^3 {+} 5x^4 {+} 3x^5
{+} x^6)\\
&=as^{-3}(r {+} s)^2(r^2{ +} s^2)(r^2{ +} rs {+} s^2) \\
&\quad- as^{-3}r^{-3}(r{ +} s)^2(r^2 {+ }s^2)(r^2{ +} rs{ +}
s^2)(r^{2}s{ +}
2r^3{ +} rs^2) \\
&\quad+ as^{-3}r^{-3}(r^2{ +} s^2)(r^2 {+ }rs{ +} s^2)^2(r{ +}
s)^3 -
ar^{-3}s^{-2}(r{ +} s)^2(r^2 {+} s^2)(r^2 {+} rs{ +} 2s^2)\\
&\quad + ar^{-3}(r {+}
s)^2(r^2 {+} s^2)(r^2 {+} rs{ +} s^2)\\
&= a(r {+} s)^2(r^2 {+} s^2)(r^2 {+} rs{ +} s^2)[s^{-3}-
s^{-3}r^{-3}(r^{2}s {+} rs^2{ +} 2r^3)\\
&\quad + r^{-3}s^{-3}(r^2 {+ }rs {+ }s^2)(r {+} s) -
r^{-3}s^{-2}(r^2{ +} rs{ +} 2s^2) +
r^{-3}]\\
&= 0.
\end{split}
\end{equation*}

(v) $X = f_{1}f_{2}f_{1}^{3}$. By calculation, we get $120$ relevant
terms of $\triangle ^{(4)}(X)$ in ($\ast$) and their pairing values
of ($\ast$) listed as in TABULAR 5 of Appendix.

If we sum up all the pairing-values in TABULAR $5$, then we get the
paring value of ($\ast$):
\begin{equation*}
\begin{split}
&ay^3(1{+} 3x {+} 5x^2 {+} 6x^3 {+} 5x^4 {+} 3x^5 {+} x^6) - a(r {+}
s)(r^2 {+} s^2)[\,\bar{x}y^{3}x^3(1 {+} 2x\\
&\quad {+} 2x^2 {+} x^3) {+} y^2(1 {+} 3x {+} 5x^2 {+} 5x^3 {+} 3x^4
{+} x^5)\,] +
ars(r^2 {+} s^2)(r^2 {+} rs {+} s^2)\cdot\\
&\qquad[\,y(1 {+} 3x {+} 4x^2 {+} 3x^3 {+} x^4) + \bar{x}y^{2}x^2(1
{+} 3x {+}
4x^2 {+} 3x^3 {+} x^4)\,]\\
&\quad - a(rs)^3(r{+}s)(r^2 {+} s^2)[\,\bar{x}yx(1
{+} 3x {+} 5x^2 {+} 5x^3 {+} 3x^4 {+} x^5) {+} 1 {+} 2x {+} 2x^2 {+} x^3\,]\\
&\quad +
ar^{6}s^6\bar{x}(1 {+} 3x {+} 5x^2 {+} 6x^3 {+} 5x^4 {+} 3x^5 {+} x^6)\\
&=a(r{+}s)^3(r^2{+}s^2)(r^2{+}rs{+}s^2) - as(r{+}s)^2(r^2{+}
s^2)(r^2{+}rs{+} s^2)(r^{2}{+}rs{+}2s^2)\\
&\quad + a(r^2{+}s^2)(r^2{+}rs {+}s^2)^2(r{+}s)^3 -
ar(r{+}s)^2(r^2{+}s^2)(r^2{+}rs{+}s^2)(2r^2{ +} rs
{+} s^2)\\
&= a(r {+} s)^2(r^2 {+} s^2)(r^2 {+} rs {+} s^2)[\,r^3 {+} s^3 {-}
s(r^{2} {+}rs {+} 2s^2)\\
&\quad {+} (r {+} s)(r^2 {+} rs {+} s^2) {-} r(2r^2 {+}
rs {+} s^2)\,]\\
&= 0.
\end{split}
\end{equation*}

Up to now, these five cases of $\Delta^{(4)}(X)$ have been checked.
The proof is completed by checking that the relations in $B'^{cop}$
are preserved for $G_2$.
\end{proof}

\begin{defi}
For any two Hopf algebras $\mathcal{A}$ and $\mathcal{U}$ connected
by a skew-dual pairing $\langle\, ,\rangle$ one may form the
Drinfel'd quantum double $\mathcal{D}(\mathcal{A},\,\mathcal{U})$ as
in $[KS, 3.2]$, which is a Hopf algebra whose underlying coalgebra
is $\mathcal{A}\otimes\mathcal{U}$ with the tensor product coalgebra
structure, whose algebra structure is defined by
$$(a\otimes
f)(a'\otimes f') = \sum\langle
\mathcal{S}_{\mathcal{U}}(f_{(1)}),\,a'_{(1)}\rangle\langle
(f_{(3)}),\,a'_{(3)}\rangle aa'_{(2)}\otimes f_{(2)}f',\leqno(6)$$
for $a, \,a'\in\mathcal{A}$ and $f,\, f'\in\mathcal{U}$, and whose
antipode $S$ is given by
$$S(a\otimes f) = (1\otimes
\mathcal{S}_{\mathcal{U}}(f))(\mathcal{S}_{\mathcal{A}}(a)\otimes
1). \leqno(7)$$
\end{defi}

Clearly, both mappings $\mathcal{A}\ni  a \mapsto a\otimes 1 \in
\mathcal{D}(\mathcal{A},\,\mathcal{U})$ and $\mathcal{U}\ni  f
\mapsto 1\otimes f \in \mathcal{D}(\mathcal{A},\,\mathcal{U})$ are
injective Hopf algebra homomorphisms. Denote the image $a\otimes 1$
 of $a$ in
$\mathcal{D}(\mathcal{A},\,\mathcal{U})$ by $\hat{a}$ and the image
$1\otimes f$ of $f$ by $\hat{f}$. By (6), we have the following
cross relations between elements $\hat{a}$ (for $a\in \mathcal{A}$)
and $\hat{f}$ (for $f\in \mathcal{U}$) in the algebra
$\mathcal{D}(\mathcal{A},\,\mathcal{U})$:
\begin{gather*}
\hat{f}\hat{a} = \sum\langle
\mathcal{S}_{\mathcal{U}}(f_{(1)}),\,a_{(1)}\rangle\langle
(f_{(3)}),\,a_{(3)}\rangle \hat{a}_{(2)}\hat{f}_{(2)},\tag{8}\\
\sum\langle f_{(1)},\,a_{(1)}\rangle \hat{f}_{(2)}\hat{a}_{(2)} =
\sum\hat{a}_{(1)}\hat{f}_{(1)}\langle
f_{(2)},\,a_{(2)}\rangle.\tag{9}
\end{gather*}
In fact, as an algebra the double
$\mathcal{D}(\mathcal{A},\,\mathcal{U})$ is the universal algebra
generated by the algebras $\mathcal{A}$ and $\mathcal{U}$ with
cross relations (8) or, equivalently, (9).

\begin{theorem}
The two-parameter quantum group $U_{r,\,s}(G_2)$ is isomorphic to
the Drinfel'd quantum double
$\mathcal{D}(\mathcal{B},\,\mathcal{B}')$.
\end{theorem}

The proof is the same as that of [BGH1, Theorem 2.5].

\begin{remark} \ The proofs of Proposition 2.3 and Theorem 2.5
show the compatibility of the defining relations of $U_{r,s}(G_2)$,
where the proof of Theorem 2.5 indicates that the cross relations
between $\mathcal{B}$ and $\mathcal{B}'$ are precisely half the ones
appearing in ($G1$)--($G4$), and the proof of Proposition 2.3 then
shows the compatibility of the remaining relations appearing in
$\mathcal{B}$ and $\mathcal{B}'$ including the other half of
($G1$)--($G4$) and the $(r,s)$-Serre relations ($G5$)--($G6$).
\end{remark}

\bigskip
\section{Lusztig's Symmetries from
$U_{r,\,s}(G_2)$ to $U_{s^{-1},\,r^{-1}}(G_2)$}
\medskip

As we did in [BGH1] for the classical types $A,\,B,\,C,\,D$, we call
$(U_{s^{-1},r^{-1}}(G_2)$, $\lg\,|\,\rg)$ the associated quantum
group corresponding to $(U_{r,s}(G_2), \lg\,,\rg)$, where the
pairing $\lg\om_i'|\,\om_j\rg$ is defined by replacing $(r, s)$ with
$(s^{-1}, r^{-1})$ in the defining formula for $\lg
\om_i',\om_j\rg$. Obviously,
$$
\lg \om_i'|\,\om_j\rg=\lg\om_j',\om_i\rg.
$$

We now study Lusztig's symmetry property between $(U_{r,s}(G_2),
\lg\,,\rg)$ and its associated object $(U_{s^{-1},r^{-1}}(G_2),
\lg\,|\,\rg)$, which indeed indicates the difference in structures
between the two-parameter quantum group introduced above and the
usual one-parameter quantum group of Drinfel'd-Jimbo type.

To define the Lusztig's symmetries, we introduce the notation of
divided-power elements (in $(U_{s^{-1},r^{-1}}(G_2),
\lg\,|\,\rg)$\,). For any nonnegative integer $k\in \Bbb N$, set
$$
\lg k\rg_i=\frac{s_i^{-k}-r_i^{-k}}{s_i^{-1}-r_i^{-1}},\qquad \lg
k\rg_i!=\lg 1\rg_i\lg 2\rg_i\cdots\lg k\rg_i,
$$
and for any element $e_i, f_i\in (U_{s^{-1},r^{-1}}(G_2),
\lg\,|\,\rg)$, define the divided-power elements
$$
 e_i^{(k)}=e_i^{k}/\lg k\rg_i!\,,\qquad f_i^{(k)}=f_i^{k}/\lg k\rg_i!\,.
$$

\begin{defi}
To every $i$ $(i= 1, 2)$, there corresponds a $\mathbb{Q}$-linear
mapping $\mathcal{T}_i: \ (U_{r,s}(G_2),\,\langle\,,\rangle)
\longrightarrow (U_{s^{-1},\,r^{-1}}(G_2),\,\langle\, |\,\rangle)$
such that $\mathcal T_i(r)=s^{-1}$, $\mathcal T_i(s)=r^{-1}$,
which acts on the generators $\omega_j,\, \omega_j',\, e_j,\,
f_j\; (1 \leq j \leq 2)$ as

\begin{gather*} \mathcal
T_i(\om_j)=\om_j\,\om_i^{-a_{ij}},\qquad \mathcal
T_i(\om_j')=\om_j'\,{\om_i'}^{-a_{ij}},\\
\mathcal T_i(e_i)=-\,{\om_i'}^{-1}f_i,\qquad \mathcal
T_i(f_i)=-(r_is_i)\,e_i\,\om_i^{-1},\\
\end{gather*}
and for $i\ne j$,
\begin{gather*} \mathcal
T_i(e_j)=\sum_{\nu=0}^{-a_{ij}}(-1)^{\nu}(rs)^{\frac{\nu}2(-a_{ij}-\nu)}\lg\om_j',\om_i\rg^{-\nu}
\lg\om_i',\om_i\rg^{\frac{\nu}{2}(1{+}a_{ij})}
e_i^{(\nu)}e_j\,e_i^{(-a_{ij}-\nu)},\\
\mathcal
T_i(f_j)=(r_js_j)^{\delta_{ij}^+}\sum_{\nu=0}^{-a_{ij}}(-1)^{\nu}(rs)^{\frac{\nu}2(-a_{ij}-\nu)}
\lg\om_i',\om_j\rg^\nu\lg\om_i',\om_i\rg^{-\frac{\nu}{2}(1{+}a_{ij})}
f_i^{(-a_{ij}-\nu)}f_j\,f_i^{(\nu)},\\
\end{gather*}
here $(a_{ij})$ is the Cartan matrix of the simple Lie algebra
$\frak g$ of type $G_2$, and for any $i\ne j$,
$$
\delta_{ij}^+=\begin{cases} 2, & \text{ if }\ i<j,\ \&\ a_{ij}\ne0, \\
1, & \text{ otherwise }.\\
\end{cases}
$$
\end{defi}

\begin{lemm}
$\mathcal{T}_i \,(i = 1, 2)$ preserves the defining relations
$(G1)$--$(G3)$ of $U_{r,s}(G_2)$ into its associated object
$U_{s^{-1},r^{-1}}(G_2)$.
\end{lemm}
\begin{proof} \ For $G_2$, we have
\begin{gather*}
\langle \omega_1',\,\omega_1\rangle = rs^{-1} = \langle \omega_1'
|\,\omega_1\rangle, \qquad\quad \langle
\omega_1',\,\omega_2\rangle
= r^{-3} = \langle \omega_2' |\,\omega_1\rangle, \\
\langle \omega_2',\,\omega_1\rangle = s^{3} = \langle \omega_1'
|\,\omega_2\rangle, \qquad\quad \langle
\omega_2',\,\omega_2\rangle = r^{3}s^{-3} = \langle \omega_2'
|\,\omega_2\rangle.
\end{gather*}

We show that $\mathcal{T}_1, \mathcal{T}_ 2$ preserve the defining
relations ($G1$)--($G3$). ($G1$) are automatically satisfied. To
check ($G2$) and ($G3$):  first of all, by direct calculation, we
have $\mathcal{T}_k(\langle\omega_i',\,\omega_j\rangle) =
\langle\mathcal{T}_k(\omega_i'),\,\mathcal{T}_k(\omega_j)\rangle =
\langle\omega_j',\,\omega_i\rangle =
\langle\omega_i'|\,\omega_j\rangle,$ for $i, j, k\in \{1,\,2\}.$
This fact ensures that $\mathcal{T}_k (k = 1, 2)$ preserve ($G2$)
and ($G3$), that is,
\begin{gather*}
\mathcal{T}_k(\omega_j)\mathcal{T}_k(e_i)\mathcal{T}_k(\omega_j)^{-1}
= \langle\omega_i'|\,\omega_j\rangle\mathcal{T}_k(e_i),\quad
\mathcal{T}_k(\omega_j)\mathcal{T}_k(f_i)\mathcal{T}_k(\omega_j)^{-1}
=\langle\omega_i'|\,\omega_j\rangle^{-1}\mathcal{T}_k(f_i),\\
\mathcal{T}_k(\omega_j')\mathcal{T}_k(e_i)\mathcal{T}_k(\omega_j')^{-1}
= \langle\omega_j'|\,\omega_i\rangle^{-1}\mathcal{T}_k(e_i),\quad
\mathcal{T}_k(\omega_j')\mathcal{T}_k(f_i)\mathcal{T}_k(\omega_j')^{-1}
= \langle\omega_j'|\,\omega_i\rangle\mathcal{T}_k(f_i),
\end{gather*}
where checking other three identities is equivalent to checking
the first one.
\end{proof}

\begin{lemm}
$\mathcal{T}_i\, (i= 1, 2)$ preserves the defining relations
$(G4)$ into its associated object $U_{s^{-1},r^{-1}}(G_2)$.
\end{lemm}
\begin{proof} \ Put $\Delta=r^2 + rs + s^2$. To check
($G4$): \ for $i= 1, 2$, we have
$$
[\mathcal{T}_i(e_i),\,\mathcal{T}_i(f_i)] =
(r_{i}s_i)\omega_i'^{-1}(f_{i}e_i - e_{i}f_i)\omega_i^{-1}
 = \mathcal{T}_i([e_i,f_i]),
$$
\begin{equation*}
\begin{split}
[\mathcal{T}_2&(e_1),\mathcal{T}_2(f_1)] = [e_{1}e_2 -
r^{3}e_{2}e_1,\,rs(f_{2}f_1 - s^{3}f_{1}f_2)]\\
&=rs\{f_2[e_1,f_1]e_2 + e_1[e_2,f_2]f_1 -
r^3([e_2,f_2]f_{1}e_1 + e_{2}f_2[e_1,f_1]) \\
&\quad- s^3([e_1,f_1]f_{2}e_2+ e_{1}f_1[e_2,f_2]) +
(rs)^3(e_2[e_1,f_1]f_2 +
f_1[e_2,\,f_2]e_1)\}\\ 
&=\frac{\omega_2\omega_1 - \omega_2'\omega_1'}{s^{-1} - r^{-1}}
=\frac{\mathcal{T}_2(\omega_1) - \mathcal{T}_2(\omega_1')}{s^{-1}
- r^{-1}}= \mathcal{T}_2([e_1,f_1]),
\end{split}
\end{equation*}
and as for
\begin{equation*}
\begin{split}
[\mathcal{T}_1(e_2),\mathcal{T}_1(f_2)]& = \frac{r^3s^3}{(r +
s)^2\Delta^2}\big[\,(rs^2)^3e_{2}e_1^{3} - rs^{3}\Delta
e_{1}e_{2}e_1^{2} + s\Delta e_1^{2}e_{2}e_{1} - e_1^{3}e_{2},\\
&\hskip2cm\quad (r^2s)^3f_1^{3}f_{2} - sr^{3}\Delta
f_1^{2}f_{2}f_{1} + r\Delta f_{1}f_{2}f_1^{2} - f_{2}f_1^{3}\,\big],
\end{split}
\end{equation*}
we have to show that the above bracket on the right-hand side is
equal to
$$
\Delta(r + s)^2\cdot\frac{\omega_2\omega_1^3 -
\omega_2'\omega_1'^3}{r - s}.
$$

To do so, we introduce the notations of ``quantum root vectors" in
terms of adjoint actions, as follows:
\begin{equation*}
\begin{split}
E_{12}&=(\text{ad}_{l}e_1)(e_2) = e_{1}e_2 -
s^3e_2e_1,\\
F_{12}&=(\text{ad}_rf_1)(f_2)=f_2f_1-r^3f_1f_2,\\
E_{112}&=(\text{ad}_{l}e_1)^2(e_2) = e_{1}E_{12} - rs^2E_{12}e_1, \\ 
F_{112}&=(\text{ad}_{r}f_1)^2(f_2) = F_{12}f_1 - r^2sf_1F_{12}, \\
E_{1112}&=(\text{ad}_{l}e_1)^3(e_2) = e_{1}^3e_2 - s\Delta
e_1^2e_2e_1 + rs^3\Delta e_1e_2e_{1}^2 - (rs^2)^3e_2e_{1}^3,\\
F_{1112}&=(\text{ad}_{r}f_1)^3(f_2) =f_{2}f_1^{3} -r\Delta
f_{1}f_{2}f_1^{2}  +sr^{3}\Delta f_1^{2}f_{2}f_{1}
-(r^2s)^3f_1^{3}f_{2}.\\
\end{split}
\end{equation*}
That is, we need to verify that
$$
[\,E_{1112},F_{1112}\,] =\Delta(r +
s)^2\cdot \frac{\omega_2\omega_1^3 - \omega_2'\omega_1'^3}{r -
s}.
$$

By direct calculation using the Leibniz rule, we have
\begin{gather*}
[e_1,F_{12}] = -\Delta\omega_1f_2,\qquad [e_2,F_{12}] =
f_1\omega_2',\\
[E_{12},f_1] = -\Delta e_2\omega_1',\qquad [E_{12},f_2] =
\omega_2e_1,\\
[E_{12},F_{12}] = \frac{\omega_1\omega_2 - \omega_1'\omega_2'}{r -
s},\\
[e_1,F_{112}] = -(r + s)^2\omega_1F_{12},\qquad [e_2,F_{112}] =
s(s^2 - r^2)f_1^2\omega_2',\\
[E_{112},f_1] = -(r + s)^2E_{12}\omega_1',\qquad [E_{112},f_2] =
r(r^2 - s^2)\omega_2e_1^2,\\
[E_{112},F_{12}] = (r + s)^2\omega_1\omega_2e_1,\qquad
[E_{12},F_{112}] = (r + s)^2f_1\omega_1'\omega_2', \\
[E_{112},F_{112}] =(r + s)^2\cdot \frac{\omega_1^2\omega_2 -
\omega_1'^2\omega_2'}{r - s},\\
\end{gather*}
as well as
\begin{equation*}
\begin{split}
[e_1&,F_{1112}]=[e_1,
F_{112}f_1-rs^2f_1F_{112}]=-\Delta\omega_1F_{112},\\
[E_{112}&,F_{1112}]=[E_{112},F_{112}f_1-rs^2f_1F_{112}]\\
&=[E_{112},F_{112}]f_1-rs^2f_1[E_{112},F_{112}]+F_{112}[E_{112},f_1]-rs^2[E_{112},f_1]F_{112}\\
&=\Delta(r+s)^2 f_1\omega_1'^2\omega_2',\\
\end{split}
\end{equation*}
\begin{equation*}
\begin{split}
[E_{1112}&,F_{1112}]=[e_1E_{112}-r^2sE_{112}e_1,F_{1112}]\\
&=[e_1,F_{1112}]E_{112}-r^2sE_{112}[e_1,F_{1112}]+e_1[E_{112},F_{1112}]-r^2s[E_{112},F_{1112}]e_1\\
&=\Delta\omega_1[E_{112},F_{112}]+\Delta(r+s)^2[e_1,f_1]\omega_1'^2\omega_2'\\
&=\Delta(r + s)^2\cdot\frac{\omega_2\omega_1^3 -
\omega_2'\omega_1'^3}{r - s}.
\end{split}
\end{equation*}

Thus, we arrive at $[\mathcal{T}_1{(e_2)},\mathcal{T}_1{(f_2)}]
 = \mathcal{T}_1([e_2,\,f_2]) \in U_{s^{-1},\,r^{-1}}(G_2)$.
\end{proof}

\begin{lemm}
$\mathcal{T}_2$ preserves the $(r,s)$-Serre relations $(G5)_1,
(G6)_1$ into its associated object $U_{s^{-1},r^{-1}}(G_2)$:
\begin{gather*}
\mathcal{T}_2(e_2)^{2}\mathcal{T}_2(e_1) - (r^{3} +
s^{3})\mathcal{T}_2(e_2)\mathcal{T}_2(e_1)\mathcal{T}_2(e_2) +
(rs)^3\mathcal{T}_2(e_1)\mathcal{T}_2(e_2)^{2} = 0,\tag{1}\\
\mathcal{T}_2(f_1)\mathcal{T}_2(f_2)^{2} - (r^{3} +
s^{3})\mathcal{T}_2(f_2)\mathcal{T}_2(f_1)\mathcal{T}_2(f_2) +
(rs)^3\mathcal{T}_2(f_1)\mathcal{T}_2(f_2)^{2} = 0.\tag{2}
\end{gather*}
\end{lemm}
\begin{proof} \ For the degree $2$ $(r,s)$-Serre
relation $(G5)_1$
$$e_2^{2}e_1 - (r^{-3} + s^{-3})e_2e_1e_2 +
r^{-3}s^{-3}e_1e_2^{2} = 0,$$ observe that
$$
\mathcal{T}_2(e_1)\mathcal{T}_2(e_2) =
r^{-3}\mathcal{T}_2(e_2)\mathcal{T}_2(e_1) -
r^{-3}e_1,\qquad\mathcal{T}_2(e_2)e_1 = s^3e_1\mathcal{T}_2(e_2).
\leqno(3)$$

Making $\mathcal{T}_2$ act algebraically on the left-hand side of
$(G5)_1$, we have
\begin{equation*}
\begin{split}
\mathcal{T}_2&(e_2)^{2}\mathcal{T}_2(e_1) - (r^{3} +
s^{3})\mathcal{T}_2(e_2)\mathcal{T}_2(e_1)\mathcal{T}_2(e_2) +
(rs)^3\mathcal{T}_2(e_1)\mathcal{T}_2(e_2)^{2} \\
&\qquad=
\mathcal{T}_2(e_2)r^3(\mathcal{T}_2(e_1)\mathcal{T}_2(e_2) +
r^{-3}e_1) - (r^{3} +
s^{3})\mathcal{T}_2(e_2)\mathcal{T}_2(e_1)\mathcal{T}_2(e_2)\\
&\qquad\quad+ (rs)^3(r^{-3}\mathcal{T}_2(e_2)\mathcal{T}_2(e_1) -
r^{-3}e_1)\mathcal{T}_2(e_2)\\
&\qquad = 0,
\end{split}
\end{equation*}
proving (1). The proof of (2) is similar.
\end{proof}

To prove that $\mathcal{T}_1$ preserves the Serre relations, we need
three auxiliary lemmas.
\begin{lemm}
In the notation in Lemma 3.3, we have
$$[E_{1112}E_{112}-r^3E_{112}E_{1112},\,f_2] = 0.$$
\end{lemm}
\begin{proof} \ Since $e_1E_{1112} - r^3E_{1112}e_1=\text{ad}_l(e_1)^4(e_2)=0$ (Serre relation), and
\begin{equation*}
\begin{split}
[E_{1112},\,f_2]& = [e_1E_{112}- r^2sE_{112}e_1,\,f_2] =
e_1[E_{112},\,f_2] - r^2s[E_{112},\,f_2]e_1\\
& = r^3(r{-}s)(r^2{-}s^2)\,\omega_2e_1^3,
\end{split}
\end{equation*}
we obtain
\begin{equation*}
\begin{split}
[E_{1112}&E_{112} - r^3E_{112}E_{1112},\,f_2] =
E_{1112}[E_{112},\,f_2] + [E_{1112},\,f_2]E_{112}\\
&\qquad - r^3\big ( E_{112}[E_{1112},\,f_2] +
[E_{112},\,f_2]E_{1112}\big )\\
&= r^3(r{-}s)(r^2{-} s^2)\,\omega_2\,\bigl(e_1^3E_{112}-r\Delta
e_1E_{1112}e_1-(r^2s)^3E_{112}e_1^3\bigr) \\
&= r^3(r{-}s)(r^2 {-} s^2)\,\omega_2\,\bigl(e_1^3E_{112}-r\Delta
e_1^2E_{112}e_1+r^3s\Delta
e_1E_{112}e_1^2-(r^2s)^3E_{112}e_1^3\bigr)\\
&= r^3(r{-}s)(r^2 {-} s^2)\,\omega_2\,\bigl(e_1\cdot({\mathcal
{SR}})-rs^2({\mathcal {SR}})\cdot e_1\bigr)\\
&=0,
\end{split}
\end{equation*}
where ${\mathcal {(SR)}}$ denotes the left-hand-side presentation
of the $(r,s)$-Serre relation $(G5)_2$
$$
e_1^2E_{112} - r^2(r + s)e_1E_{112}e_1 + r^5sE_{112}e_1^2=0,
$$
and we used the replacement $E_{1112}=e_1E_{112}-r^2sE_{112}e_1$
in the third equality.
\end{proof}

\begin{lemm}
In the notation in Lemma 3.3, we have
$$[E_{1112}E_{112} - r^3E_{112}E_{1112},\,f_1] = 0.$$
\end{lemm}
\begin{proof} \ It is easy to check that $[E_{1112}, f_1]=-\Delta
E_{112}\omega_1'$. Thus
\begin{equation*}
\begin{split}
[E_{1112}&E_{112} - r^3E_{112}E_{1112},\,f_1] =
E_{1112}[E_{112},\,f_1] + [E_{1112},\,f_1]E_{112}\\
&\qquad - r^3\big ( E_{112}[E_{1112},\,f_1] +
[E_{112},\,f_1]E_{1112}\big )\\
&=(r+s)\,\bigl[(r + s)((rs)^3E_{12}E_{1112} - E_{1112}E_{12}) +
r(r - s)\Delta E_{112}^2\bigr]\omega_1'.
\end{split}
\end{equation*}

It suffices to show that
$$
E_{1112}E_{12} = (rs)^3E_{12}E_{1112} + r(r - s)(r + s)^{-1}\Delta
E_{112}^2.\leqno(4) $$

At first, we note that the $(r,s)$-Serre relation $(G5)_1$ is
equivalent to
$$E_{12}e_2 = r^3e_2E_{12}.$$
As $e_1e_2=E_{12}+s^3e_2e_1$, we get
\begin{equation*}
\begin{split}
E_{112}e_2
&=(e_1E_{12}-rs^2E_{12}e_1)e_2=r^3e_1e_2E_{12}-rs^2E_{12}e_1e_2\\
&=r^3(E_{12}+s^3e_2e_1)E_{12}-rs^2E_{12}(E_{12}+s^3e_2e_1)\\
&=r(r^2{-}s^2)E_{12}^2+(rs)^3e_2(e_1E_{12}-rs^2E_{12}e_1)\\
&=r(r^2{-}s^2)E_{12}^2+(rs)^3e_2E_{112}.
\end{split}
\end{equation*}

Next, we claim
$$E_{1112}e_2 = (rs^2)^3e_2E_{1112} - r(rs {-} r^2 {+}
s^2)E_{112}E_{12} + (rs)^2(r^2 {+} rs {-} s^2)E_{12}E_{112}.$$

Indeed, since $E_{1112} = e_1E_{112} - r^2sE_{112}e_1, E_{112} =
e_1E_{12} - rs^2E_{12}e_1$, $e_1e_2=E_{12}+ s^3e_2e_1$, we have
\begin{equation*}
\begin{split}
E_{1112}e_2&=e_1(E_{112}e_2)-r^2sE_{112}(e_1e_2)\\
&=r(r^2-s^2)e_1E_{12}^2+(rs)^3(e_1e_2)E_{112}-r^2sE_{112}(e_1e_2)\\
&=r(r^2-s^2)e_1E_{12}^2+(rs)^3E_{12}E_{112}+(rs^2)^3e_2e_1E_{112}-r^2sE_{112}E_{12}\\
&\qquad-(rs^2)^2(E_{112}e_2)e_1\\
&=r(r^2-s^2)E_{112}E_{12}+(rs)^2(r^2-s^2)E_{12}e_1E_{12}+(rs)^3E_{12}E_{112}\\
&\qquad+(rs^2)^3e_2e_1E_{112}-r^2sE_{112}E_{12}-r^3s^4(r^2-s^2)E_{12}^2e_1-r^5s^7e_2E_{112}e_1\\
&=(rs^2)^3e_2E_{1112} - r(rs {-} r^2 {+} s^2)E_{112}E_{12} +
(rs)^2(r^2 {+} rs {-} s^2)E_{12}E_{112}.
\end{split}
\end{equation*}

To prove (4), we first note that
\begin{equation*}
\begin{split}
[(r +& s)((rs)^3E_{12}E_{1112} - E_{1112}E_{12}) + r(r -
s)\Delta E_{112}^2,\,f_1]\\
&=(r + s)(rs)^3\bigl(E_{12}[E_{1112},\,f_1] +
[E_{12},\,f_1]E_{1112}\bigr)\\
&\quad - (r + s)\bigl(E_{1112}[E_{12},\,f_1] +
[E_{1112},\,f_1]E_{12}\bigr) \\
&\quad+ r(r - s)\Delta(E_{112}[E_{112},f_1]+[E_{112},\,f_1]E_{112})\\
&=-(r + s)(rs)^3\Delta\bigl(E_{12}E_{112} +
s^3e_2E_{1112}\bigr)\omega_1'\\
&\quad+(r+s)\Delta\bigl(E_{1112}e_2+r^2sE_{112}E_{12}\bigl)\omega_1'\\
&\quad-r(r-s)(r+s)^2\Delta\bigl(E_{112}E_{12}+rs^2E_{12}E_{112}\bigr)\omega_1',
\end{split}
\end{equation*}
which vanishes by the preceding identity. Second, instead of $f_1$
by $f_2$ in the above formula, we get
\begin{equation*}
\begin{split}
[(r +& s)((rs)^3E_{12}E_{1112} - E_{1112}E_{12}) + r(r -
s)\Delta E_{112}^2,\,f_2]\\
&=(r + s)(rs)^3\bigl(r^3(r^2-s^2)(r-s)E_{12}\omega_2e_1^3+\omega_2e_1E_{1112}\bigr)\\
&\quad - (r + s)\bigl(E_{1112}\omega_2e_1
+ r^3(r^2-s^2)(r-s)\omega_2e_1^3E_{12}\bigr)\\
&\quad + r^2(r - s)(r^2-s^2)\Delta\bigl(E_{112}\omega_2e_1^2+\omega_2e_1^2E_{112}\bigr)\\
&=(r +
s)(rs)^3\omega_2\bigl((rs)^3(r^2-s^2)(r-s)E_{12}e_1^3+e_1E_{1112}\bigr)\\
&\quad-(r +
s)\omega_2\bigl((r^2s)^3E_{1112}e_1+r^3(r^2-s^2)(r-s)e_1^3E_{12}\bigr)\\
&\quad+r^2(r -
s)(r^2-s^2)\Delta\omega_2\bigl((rs)^3E_{112}e_1^2+e_1^2E_{112}\bigr)\\
&=r^2(r -
s)(r^2-s^2)\Delta\omega_2\bigl((rs)^3E_{112}e_1^2+e_1^2E_{112}\bigr)\\
&\quad-r^3(r^2-s^2)^2\omega_2\bigl(e_1^3E_{12}-(rs^2)^3E_{12}e_1^3\bigr)\\
&=r^2(r -
s)(r^2-s^2)\Delta\omega_2\bigl((rs)^3E_{112}e_1^2+e_1^2E_{112}\bigr)\\
&\quad-r^3(r^2-s^2)^2\omega_2\bigl(e_1^2E_{112}+rs^2e_1E_{112}e_1+(rs^2)^2E_{112}e_1^2\bigr)\\
&=(rs)^2(r -
s)(r^2-s^2)\omega_2\bigl(e_1^2E_{112}-r^2(r+s)e_1E_{112}e_1+r^5sE_{112}e_1^2\bigr)\\
&=(rs)^2(r - s)(r^2-s^2)\omega_2(\text{ad}_le_1)^2(E_{112})\\
&=(rs)^2(r - s)(r^2-s^2)\omega_2(\text{ad}_le_1)^4(e_2)\\
&=0.
\end{split}
\end{equation*}
Then, through an argument similar to the one used in the deduction
of [BKL, Lemma 3.4], we get (4).
\end{proof}

By [BKL, Lemma 3.4], Lemmas 3.5 and 3.6 imply:
\begin{lemm}
$$
E_{1112}E_{112}-r^3E_{112}E_{1112}=0.$$
\end{lemm}

\begin{lemm}
$\mathcal{T}_1$ preserves the $(r,s)$-Serre relations $(G5)_1,
(G6)_1$ into its associated object $U_{s^{-1},r^{-1}}(G_2)$, i.e.,
\begin{gather*}
\mathcal{T}_1(e_2)^{2}\mathcal{T}_1(e_1) - (r^{3} +
s^{3})\mathcal{T}_1(e_2)\mathcal{T}_1(e_1)\mathcal{T}_1(e_2) +
(rs)^3\mathcal{T}_1(e_1)\mathcal{T}_1(e_2)^{2} = 0,\tag{5}\\
\mathcal{T}_1(f_1)\mathcal{T}_1(f_2)^{2} - (r^{3} +
s^{3})\mathcal{T}_1(f_2)\mathcal{T}_1(f_1)\mathcal{T}_1(f_2) +
(rs)^3\mathcal{T}_1(f_2)^{2}\mathcal{T}_1(f_1) = 0.\tag{6}
\end{gather*}
\end{lemm}
\begin{proof} \ By direct calculation, we have
\begin{equation*}
\begin{split}
\mathcal{T}_1(e_2)\mathcal{T}_1(e_1) &= \Big[-\frac{1}{s^3(r +
s)\Delta}E_{1112}\Big]\cdot(-\omega_1'^{-1}f_1)\\
&= s^3\mathcal{T}_1(e_1)\mathcal{T}_1(e_2) - \frac{1}{rs^2(r +
s)}E_{112}.
\end{split}\tag{7}
\end{equation*}
Hence, to prove (5) is equivalent to prove
$$
\mathcal T_1(e_2)E_{112}-r^3E_{112}\mathcal T_1(e_2)=0.$$ However,
the latter is given by Lemma 3.7.

The proof of (6) is analogous.
\end{proof}

To prove that $\mathcal{T}_2$ preserves the Serre relations, we also
need auxiliary lemmas. Write
$$E_{21}:= (ad_{l}e_2)(e_1) = e_2e_{1} - r^{-3}e_1e_2,$$
and note that $(G5)_1$ is equivalent to $(ad_{l}e_2)(E_{21}) =
e_2E_{21} - s^{-3}E_{21}e_2 = 0$, i.e., $E_{21}e_2=s^3e_2E_{21}$.

\begin{lemm}
$$\bigl[\,e_1E_{21}^3 - s\Delta E_{21}e_1E_{21}^2 + rs^3\Delta
E_{21}^2e_1E_{21} - (rs^2)^3E_{21}^3e_1,\,f_1\,\bigr] = 0.$$
\end{lemm}
\begin{proof} \
Since
\begin{gather*}
[E_{21},\,f_1] = r^{-3}\Delta e_2\omega_1,\qquad\qquad\quad \omega_1E_{21}=rs^2E_{21}\omega_1,\\
[E_{21}^2,\,f_1] = r^{-3}s^{-1}(r + s)\Delta E_{21}e_2\omega_1,\qquad \omega_1'E_{21}=r^2sE_{21}\omega_1',\tag{8}\\
[E_{21}^3,\,f_1] = r^{-3}s^{-2}\Delta^2 E_{21}^2e_2\omega_1,
\end{gather*}
we get
\begin{equation*}
\begin{split}
\hbox{$\sum_1$}&=\frac{\omega_1 {-} \omega_1' }{r {-} s}E_{21}^3-
(rs^2)^3E_{21}^3\frac{\omega_1 {-} \omega_1' }{r {-} s} - s\Delta
E_{21}\frac{\omega_1 {-} \omega_1' }{r {-} s}E_{21}^2+rs^3 \Delta
E_{21}^2\frac{\omega_1 {-} \omega_1' }{r {-} s}E_{21} \\
&=-(rs)^3\Delta E_{21}^3\omega_1'+rs^2\Delta E^2_{21}\omega_1'E_{21}\\
&=0,
\end{split}
\end{equation*}
and
\begin{equation*}
\begin{split}
\bigl[\,e_1E_{21}^3& - s\Delta E_{21}e_1E_{21}^2 + rs^3\Delta
E_{21}^2e_1E_{21} - (rs^2)^3E_{21}^3e_1,\,f_1\,\bigr] \\
&= e_1[E_{21}^3,\,f_1] + [e_1,\,f_1]E_{21}^3 -
(rs^2)^3\bigl(E_{21}^3[e_1,\,f_1] + [E_{21}^3,\,f_1]e_1\bigr)\\
&\quad- s\Delta\bigl(E_{21}e_1[E_{21}^2,\,f_1] +
E_{21}[e_1,\,f_1]E_{21}^2 + [E_{21},\,f_1]e_1E_{21}^2\bigr)\\
&\quad+ rs^3\Delta\bigl(E_{21}^2e_1[E_{21},\,f_1] +
E_{21}^2[e_1,\,f_1]E_{21} + [E_{21}^2,\,f_1]e_1E_{21}\bigr)\\
&=r^{-3}s^{-2}\Delta^2e_1E_{21}^2e_2\omega_1 + \frac{\omega_1
{-} \omega_1' }{r {-} s}E_{21}^3\\
&\quad - s\Delta\bigl[r^{-3}s^{-1}(r {+} s)\Delta
E_{21}e_1E_{21}e_2\omega_1 + E_{21}\frac{\omega_1 {-} \omega_1'
}{r{-} s}E_{21}^2 + r^{-3}\Delta e_2\omega_1e_1E_{21}^2\bigr]\\
&\quad + rs^3\Delta\bigl[r^{-3}\Delta E_{21}^2e_1e_2\omega_1 +
E_{21}^2\frac{\omega_1 {-} \omega_1' }{r {-} s}E_{21} +
r^{-3}s^{-1}(r
{+} s)\Delta E_{21}e_2\omega_1e_1E_{21}\bigr]\\
&\quad - (rs^2)^3\bigl(E_{21}^3\frac{\omega_1 {-} \omega_1' }{r
{-}
s} + r^{-3}s^{-2}\Delta^2 E_{21}^2e_2\omega_1e_1\bigr)\\
&=\hbox{$\sum_1$}+(r^{-3}s^{-2}\Delta^2)\,\hbox{$\sum_2$}\,\omega_1=(r^{-3}s^{-2}\Delta^2)\,\hbox{$\sum_2$}\,\omega_1,
\end{split}
\end{equation*}
where
\begin{equation*}
\begin{split}
\hbox{$\sum_2$}&=e_1E_{21}^2e_2 - s^2(r + s)E_{21}e_1E_{21}e_2 -
(rs^2)^3e_2e_1E_{21}^2 \\
&\quad+ rs^5E_{21}^2e_1e_2 + r^3s^5(r + s)E_{21}e_2e_1E_{21} -
r^4s^5E_{21}^2e_2e_1.
\end{split}
\end{equation*}

We next show $\sum_2=0$. As $E_{21}e_2=s^3e_2E_{21}$ and
$e_2e_1-r^{-3}e_1e_2=E_{21}$, we get
\begin{equation*}
\begin{split}
\hbox{$\sum_2$}&=\bigl(e_1E_{21}^2e_2 -
(rs^2)^3e_2e_1E_{21}^2\bigr)+ \bigl(rs^5E_{21}^2e_1e_2 -
r^4s^5E_{21}^2e_2e_1\bigr)\\
&\quad -s^2(r + s)E_{21}e_1E_{21}e_2
+ r^3s^5(r +
s)E_{21}e_2e_1E_{21} \\
&=-(rs^2)^3E_{21}^3-r^4s^5E_{21}^3+ r^3s^5(r + s) E_{21}^3\\
&=0.
\end{split}
\end{equation*}
This completes the proof.
\end{proof}

\begin{lemm}
$$\bigl[\,e_1E_{21}^3 - s\Delta E_{21}e_1E_{21}^2 + rs^3\Delta
E_{21}^2e_1E_{21} - (rs^2)^3E_{21}^3e_1,\,f_2\,\bigr] = 0.$$
\end{lemm}
\begin{proof} \ Noting that
\begin{gather*}
[E_{21},\,f_2] = -r^{-3}\omega_2'e_1,\qquad\quad E_{21}\omega_2'=r^3\omega_2'E_{21},\\
[E_{21}^2,\,f_2] = -r^{-3}\omega_2'(e_1E_{21}+r^3E_{21}e_1),\tag{9}\\
[E_{21}^3,\,f_2] = -r^{-3}\omega_2'(e_1E_{21}^2 +
r^3E_{21}e_1E_{21} + r^6E_{21}^2e_1),
\end{gather*}
we obtain
\begin{equation*}
\begin{split}
\bigl[\,e_1E_{21}^3 &- s\Delta E_{21}e_1E_{21}^2 + rs^3\Delta
E_{21}^2e_1E_{21} - (rs^2)^3E_{21}^3e_1,\,f_2\,\bigr] \\
&= e_1[E_{21}^3,\,f_2] - s\Delta\big(E_{21}e_1[E_{21}^2,\,f_2] +
[E_{21},\,f_2]e_1E_{21}^2\big)\\
&\quad+ rs^3\Delta\big(E_{21}^2e_1[E_{21},\,f_2] +
[E_{21}^2,\,f_2]e_1E_{21}\big) - (rs^2)^3[E_{21}^3,\,f_2]e_1\\
&=-r^{-3}\omega_2'\left\{s^3e_1\bigl(e_1E_{21}^2+r^3E_{21}e_1E_{21}+r^6E_{21}^2e_1
\bigr)\right.\\
&\quad-s\Delta\bigl[(rs)^3
E_{21}e_1\bigl(e_1E_{21}+r^3E_{21}e_1\bigr)+e_1^2E_{21}^2\bigr]\\
&\quad+rs^3\Delta\bigl[(r^2s)^3E_{21}^2e_1^2+(e_1E_{21}+r^3E_{21}e_1)e_1E_{21}\bigr]\\
&\quad\left.- (rs^2)^3(e_1E_{21}^2 + r^3E_{21}e_1E_{21} +
r^6E_{21}^2e_1)e_1\right\}\\
&=-r^{-2}s\,\omega_2'\,S.
\end{split}
\end{equation*}
where
\begin{equation*}
\begin{split}
S&=(rs)^2(r^3{-}s^3)\bigl(e_1E_{21}^2e_1+E_{21}e_1^2E_{21}\bigr) +
s^2(2r^2{+} rs {+} s^2)(e_1E_{21})^2\\
& \quad- r^5s^3(2s^2 {+} rs {+} r^2)(E_{21}e_1)^2 - (r {+}
s)\bigl(e_1^2E_{21}^2 -(rs)^6E_{21}^2e_1^2\bigr).
\end{split}
\end{equation*}

It remains to prove $S=0$, which by [BKL, Lemma 3.4] is equivalent
to showing that $[S, f_1]=0=[S, f_2]$. To this end, we first
observe:
\begin{lemm}$$e_1^3E_{21}-s\Delta e_1^2E_{21}e_1+rs^3\Delta e_1E_{21}e_1^2
-(rs^2)^3E_{21}e_1^3 = 0.$$
\end{lemm}
\begin{proof} \ It is easy to see that
$$
e_1^3E_{21}-s\Delta e_1^2E_{21}e_1+rs^3\Delta e_1E_{21}e_1^2
-(rs^2)^3E_{21}e_1^3=r^{-3}(\text{ad}_le_1)^4(e_2),
$$
which is in fact the $(r,s)$-Serre relation $(G5)_2$ up to a
factor $r^{-3}$.
\end{proof}

Now set $S_i:=[S, f_i]$ for $i=1, 2$. Using (9), we obtain
\begin{equation*}
\begin{split}
S_2&=(rs)^2(r^3{-}s^3)\bigl(e_1[E_{21}^2,f_2]e_1+[E_{21},f_2]e_1^2E_{21}+E_{21}e_1^2[E_{21},f_2]\bigl)\\
&\quad+ s^2(2r^2{+} rs {+}
s^2)\bigl(e_1[E_{21},f_2]e_1E_{21}+e_1E_{21}e_1[E_{21},f_2]\bigr)\\
&\quad- r^5s^3(2s^2 {+} rs {+}
r^2)\bigl([E_{21},f_2]e_1E_{21}e_1+E_{21}e_1[E_{21},f_2]e_1\bigr)\\
&\quad-(r+s)\bigl(e_1^2[E_{21}^2,f_2]-(rs)^6[E_{21}^2,f_2]e_1^2\bigr)\\
&=(-r^{-3})\left\{(rs)^2(r^3{-}s^3)\bigl[e_1\om_2'(e_1E_{21}+r^3E_{21}e_1)e_1+
\om_2'e_1^3E_{21}+E_{21}e_1^2\om_2'e_1\bigr]\right.\\
&\quad+s^2(2r^2{+} rs {+}
s^2)\bigl(e_1\om_2'e_1^2E_{21}+e_1E_{21}e_1\om_2'e_1\bigr)\\
&\quad-r^5s^3(2s^2 {+} rs {+}
r^2)\bigl(\om_2'e_1^2E_{21}e_1+E_{21}e_1\om_2'e_1^2\bigr)\\
&\quad-\left.(r+s)\bigl[e_1^2\om_2'(e_1E_{21}+r^3E_{21}e_1)-(rs)^6\om_2'(e_1E_{21}+r^3E_{21}e_1)e_1^2\bigr]\right\}\\
&=(-r^{-3})\,\om_2'\left\{(rs)^2(r^3{-}s^3)\bigl[s^3e_1(e_1E_{21}+r^3E_{21}e_1)e_1+
e_1^3E_{21}+(rs^2)^3E_{21}e_1^3\bigr]\right.\\
&\quad+s^2(2r^2{+} rs {+}
s^2)\bigl(s^3e_1^3E_{21}+(rs^2)^3e_1E_{21}e_1^2\bigr)\\
&\quad-r^5s^3(2s^2 {+} rs {+}
r^2)\bigl(e_1^2E_{21}e_1+(rs)^3E_{21}e_1^3\bigr)\\
&\quad-\left.(r+s)\bigl[s^6e_1^2(e_1E_{21}+r^3E_{21}e_1)-(rs)^6(e_1E_{21}+r^3E_{21}e_1)e_1^2\bigr]\right\}\\
&=(-r^{-3})(rs)^2(r^3{+}s^3)\,\om_2'\left[\,e_1^3E_{21}-s\Delta
e_1^2E_{21}e_1+rs^3\Delta e_1E_{21}e_1^2
-(rs^2)^3E_{21}e_1^3\,\right]\\
&=0. \qquad\text{(by Lemma 3.11)}
\end{split}
\end{equation*}

Next we prove that $S_1=0$. Using (8) and noting that $[e_1^2,\,f_1]
= \frac{r + s}{rs}\cdot\frac{s\omega_1 - r\omega_1'}{r - s}e_1$, we
can get
\begin{equation*}
\begin{split}
S_1&=(rs)^2(r^3{-}s^3)\bigl([e_1,f_1]E_{21}^2e_1+e_1[E_{21}^2,f_1]e_1+e_1E_{21}^2[e_1,f_1]\\
&\hskip2.3cm+[E_{21},f_1]e_1^2E_{21}+E_{21}[e_1^2,f_1]E_{21}+E_{21}e_1^2[E_{21},f_1]\bigl)\\
&\quad+ s^2(2r^2{+} rs {+}
s^2)\bigl([e_1,f_1]E_{21}e_1E_{21}+e_1[E_{21},f_1]e_1E_{21}\\
&\hskip3cm+e_1E_{21}[e_1,f_1]E_{21}+e_1E_{21}e_1[E_{21},f_1]\bigr)\\
&\quad- r^5s^3(2s^2 {+} rs {+}
r^2)\bigl([E_{21},f_1]e_1E_{21}e_1+E_{21}[e_1,f_1]E_{21}e_1\\
&\hskip3cm+E_{21}e_1[E_{21},f_1]e_1+E_{21}e_1E_{21}[e_1,f_1]\bigr)\\
&\quad-(r+s)\bigl([e_1^2,f_1]E_{21}^2+e_1^2[E_{21}^2,f_1]-(rs)^6[E_{21}^2,f_1]e_1^2-(rs)^6E_{21}^2[e_1^2,f_1]\bigr)\\
&=(rs)^2(r^3{-}s^3)\left(\frac{\om_1{-}\om_1'}{r{-}s}E_{21}^2e_1+r^{-3}s^{-1}(r{+}s)\Delta
e_1E_{21}e_2\om_1e_1
+e_1E_{21}^2\frac{\om_1{-}\om_1'}{r{-}s}\right.\\
&\hskip1.8cm+\left.r^{-3}\Delta
e_2\om_1e_1^2E_{21}+\frac{r{+}s}{rs}E_{21}\frac{s\om_1{-}r\om_1'}{r{-}s}e_1E_{21}
+r^{-3}\Delta E_{21}e_1^2e_2\om_1\right)\\
&\quad+ s^2(2r^2{+} rs {+}
s^2)\left(\frac{\om_1{-}\om_1'}{r{-}s}E_{21}e_1E_{21}+r^{-3}\Delta e_1e_2\om_1e_1E_{21}\right.\\
&\hskip3cm+\left.e_1E_{21}\frac{\om_1{-}\om_1'}{r{-}s}E_{21}+r^{-3}\Delta e_1E_{21}e_1e_2\om_1\right)\\
&\quad- r^5s^3(2s^2 {+} rs {+}
r^2)\left(r^{-3}\Delta e_2\om_1e_1E_{21}e_1+E_{21}\frac{\om_1{-}\om_1'}{r{-}s}E_{21}e_1\right.\\
&\hskip3cm+\left.r^{-3}\Delta E_{21}e_1e_2\om_1e_1+E_{21}e_1E_{21}\frac{\om_1{-}\om_1'}{r{-}s}\right)\\
\end{split}
\end{equation*}
\begin{equation*}
\begin{split}
&\quad-(r{+}s)\left(\frac{r{+}s}{rs}\frac{s\om_1{-}r\om_1'}{r{-}s}e_1E_{21}^2+r^{-3}s^{-1}(r{+}s)\Delta
e_1^2E_{21}e_2\om_1 \right.\\
&\hskip1.5cm\left.-r^3s^5(r{+}s)\Delta E_{21}e_2\om_1e_1^2-(rs)^5(r{+}s)E_{21}^2\frac{s\om_1{-}r\om_1'}{r{-}s}e_1\right)\\
&=A+B+C+D,
\end{split}
\end{equation*}
where $A, B, C, D$ are given as follows.

\smallskip
Noting that $\om_1E_{21}=rs^2E_{21}\om_1$ and
$\om_1'E_{21}=r^2sE_{21}\om_1'$,  we have
\begin{equation*}
\begin{split}
A:&=(rs)^2\Delta(\om_1{-}\om_1')E_{21}^2e_1- r^5s^3(2s^2 {+} rs
{+}
r^2)E_{21}\frac{\om_1{-}\om_1'}{r{-}s}E_{21}e_1\\
&\quad+(rs)^5\frac{(r{+}s)^2}{r{-}s}E_{21}^2(s\om_1{-}r\om_1')e_1\\
&=-r^5s^4(r^3{-}s^3)E_{21}^2e_1\om_1,\\
B:&=(rs)(r{+}s)\Delta
E_{21}(s\om_1{-}r\om_1')e_1E_{21}+s^2(2r^2{+}rs{+}s^2)\frac{\om_1{-}\om_1'}{r{-}s}E_{21}e_1E_{21}\\
&\quad- r^5s^3(2s^2 {+} rs {+}
r^2)E_{21}e_1E_{21}\frac{\om_1{-}\om_1'}{r{-}s}=0,\\
C:&=(rs)^2\Delta e_1E_{21}^2(\om_1{-}\om_1')+s^2(2r^2{+} rs {+}
s^2)e_1E_{21}\frac{\om_1{-}\om_1'}{r{-}s}E_{21}\\
&\quad-\frac{(r{+}s)^2}{rs}\frac{s\om_1{-}r\om_1'}{r{-}s}e_1E_{21}^2\\
&=rs^2(r^3{-}s^3)e_1E_{21}^2\om_1,
\end{split}
\end{equation*}
furthermore, using $E_{21}e_2=s^3e_2E_{21}$,
$r^{-3}e_1e_2=e_2e_1-E_{21}$ and $e_2e_1=E_{21}+r^{-3}e_1e_2$, we
get
\begin{equation*}
\begin{split}
D:&=r^{-3}\Delta\left\{(rs)^2(r^3{-}s^3)\bigl[s^{-1}(r{+}s)e_1E_{21}e_2\om_1e_1+e_2\om_1e_1^2E_{21}+
E_{21}e_1^2e_2\om_1\bigr]\right.\\
&\quad+s^2(2r^2{+}rs{+}s^2)\bigl[e_1e_2\om_1e_1E_{21}+e_1E_{21}e_1e_2\om_1\bigr]\\
&\quad-r^5s^3(2s^2{+}rs{+}r^2)\bigl[e_2\om_1e_1E_{21}e_1+E_{21}e_1e_2\om_1e_1\bigr]\\
&\quad\left.-(r{+}s)^2s^{-1}\bigl[e_1^2E_{21}e_2\om_1-(rs)^6E_{21}e_2\om_1e_1^2\bigr]\right\}\\
&=\Delta\left\{(rs)^2(r^3{-}s^3)\bigl[(rs)^{-2}(r{+}s)e_1E_{21}e_2e_1+(e_2e_1)e_1E_{21}+
E_{21}e_1(r^{-3}e_1e_2)\bigr]\right.\\
&\quad+s^2(2r^2{+}rs{+}s^2)\bigl[r^{-1}se_1e_2e_1E_{21}+e_1E_{21}(r^{-3}e_1e_2)\bigr]\\
&\quad-r^5s^3(2s^2{+}rs{+}r^2)\bigl[(e_2e_1)E_{21}e_1+r^{-2}s^{-1}E_{21}e_1e_2e_1\bigr]\\
&\quad\left.-(r{+}s)^2s^2\bigl[e_1(r^{-3}e_1e_2)E_{21}-r^5sE_{21}(e_2e_1)e_1\bigr]\right\}\om_1\\
&=\Delta\left\{(rs)^2(r^3{-}s^3)\bigl[(rs)^{-2}(r{+}s)e_1E_{21}e_2e_1+r^{-3}e_1e_2e_1E_{21}+
E_{21}e_1e_2e_1\bigr]\right.\\
&\quad+s^2(2r^2{+}rs{+}s^2)\bigl[r^{-1}se_1e_2e_1E_{21}+e_1E_{21}e_2e_1-e_1E_{21}^2\bigr]\\
&\quad-r^5s^3(2s^2{+}rs{+}r^2)\bigl[E_{21}^2e_1+(rs)^{-3}e_1E_{21}e_2e_1+r^{-2}s^{-1}E_{21}e_1e_2e_1\bigr]\\
&\quad\left.-(r{+}s)^2s^2\bigl[e_1e_2e_1E_{21}-e_1E_{21}^2-r^5sE_{21}^2e_1 -r^2sE_{21}e_1e_2e_1\bigr]\right\}\om_1\\
&=(r^3{-}s^3)\left(r^5s^4E_{21}^2e_1-rs^2e_1E_{21}^2\right)\om_1.
\end{split}
\end{equation*}

Thus, we show $S_1=A+B+C+D=0$. This completes the proof of Lemma
3.10.
\end{proof}

The next identity is a consequence of Lemmas 3.9 and 3.10 and [BKL,
Lemma 3.4].
\begin{lemm}
$$
e_1E_{21}^3 - s\Delta E_{21}e_1 E_{21}^2 + rs^3\Delta E_{21}^2e_1
E_{21} - (rs^2)^3E_{21}^3e_1 = 0.
$$
\end{lemm}

\begin{lemm}
$\mathcal{T}_2$ preserves the $(r,s)$-Serre relations $(G5)_2,
(G6)_2$ into its associated object $U_{s^{-1},r^{-1}}(G_2)$.
\end{lemm}
\begin{proof} \ For the fourth degree $(r, s)$-Serre relation
$(G5)_2$, we have to prove that
\begin{gather*}
(rs)^6\mathcal{T}_2(e_1)^{4}\mathcal{T}_2(e_2) - (rs)^3(r {+}
s)(r^{2} {+}
s^{2})\mathcal{T}_2(e_1)^{3}\mathcal{T}_2(e_{2})\mathcal{T}_2(e_1)\\
 +\, (rs)(r^{2} {+} s^{2})(r^{2} {+} rs {+}
s^{2})\mathcal{T}_2(e_{1})^{2}\mathcal{T}_2(e_{2})\mathcal{T}_2(e_1)^2\\
-\, (r {+} s)(r^{2} {+}
s^{2})\mathcal{T}_2(e_{1})\mathcal{T}_2(e_{2})\mathcal{T}_2(e_1)^3 +
\, \mathcal{T}_2(e_{2})\mathcal{T}_2(e_1)^4
\end{gather*}
vanishes. By virtue of the commutative relation in (3), this is
equivalent to
\begin{gather*}
e_1\mathcal T_2(e_1)^3-s\Delta\mathcal T_2(e_1)e_1 \mathcal
T_2(e_1)^2+rs^3\Delta\mathcal T_2(e_1)^2e_1\mathcal
T_2(e_1)-(rs^2)^3\mathcal T_2(e_1)^3e_1=0.
\end{gather*}
However, as $\mathcal T_2(e_1)=e_1e_2-r^3e_2e_1=(-r^3)E_{21}$, the
above identity is exactly the one given by Lemma 3.12.

Similarly, we can verify that $\mathcal T_2$ preserves the
$(r,s)$-Serre relation $(G6)_2$ into its associated object
$U_{s^{-1},r^{-1}}(G_2)$.
\end{proof}

\begin{lemm}
$\mathcal{T}_1$ preserves the $(r,s)$-Serre relations $ (G5)_2,
(G6)_2$ into its associated object $U_{s^{-1},r^{-1}}(G_2)$.
\end{lemm}
\begin{proof} \ For the fourth degree $(r, s)$-Serre relation
$(G5)_2$, we have to prove that
\begin{gather*}
(rs)^6\mathcal{T}_1(e_1)^{4}\mathcal{T}_1(e_2) - (rs)^3(r {+}
s)(r^{2} {+}
s^{2})\mathcal{T}_1(e_1)^{3}\mathcal{T}_1(e_{2})\mathcal{T}_1(e_1)\\
 +\, (rs)(r^{2} {+} s^{2})(r^{2} {+} rs {+}
s^{2})\mathcal{T}_1(e_{1})^{2}\mathcal{T}_1(e_{2})\mathcal{T}_1(e_1)^2\\
-\, (r {+} s)(r^{2} {+}
s^{2})\mathcal{T}_1(e_{1})\mathcal{T}_1(e_{2})\mathcal{T}_1(e_1)^3
+ \, \mathcal{T}_1(e_{2})\mathcal{T}_1(e_1)^4 = 0.
\end{gather*}

In view of the commutation relation in (7), this is equivalent to
$$
E_{112}\mathcal{T}_1(e_1)^3 - r\Delta
\mathcal{T}_1(e_1)E_{112}\mathcal{T}_1(e_1)^2 +
r^3s\Delta\mathcal{T}_1(e_1)^2E_{112}\mathcal{T}_1(e_1) -
(r^2s)^3\mathcal{T}_1(e_1)^3E_{112}= 0.\leqno(10)
$$
We can further reduce this condition to
$$
E_{12}\mathcal{T}_1(e_1)^2 - r^2(r + s)
\mathcal{T}_1(e_1)E_{12}\mathcal{T}_1(e_1) +
r^5s\mathcal{T}_1(e_1)^2E_{12}= 0,\leqno(11)
$$
as a consequence of the commutative relation
$$
E_{112}\mathcal{T}_1(e_1) = rs^2\mathcal{T}_1(e_1)E_{112} +
r^{-1}s(r + s)^2E_{12},$$ itself arising from the equalities
$[E_{112}\,f_1] = -(r + s)^2E_{12}\omega_1'$ and
$\om_1'E_{112}=rs^2E_{112}\om_1'$.

Again, since $[E_{12},f_1] = - \Delta e_2\omega_1'$, we have
$$E_{12}\mathcal{T}_1(e_1) = r^2s\mathcal{T}_1(e_1)E_{12} + r^{-1}s\Delta
e_2,$$ by which (11) is finally reduced to $e_2\mathcal{T}_1(e_1) =
r^3 \mathcal{T}_1(e_1)e_2$, since $\mathcal
T_1(e_1)=-\om_1'^{-1}f_1$.

The proof of the second part is similar.
\end{proof}

\begin{theorem}
$\mathcal{T}_1$ and $\mathcal{T}_2$ are the Lusztig's symmetries
from $U_{r,s}(G_2)$ to its associated quantum group
$U_{s^{-1},r^{-1}}(G_2)$ as $\Bbb Q$-isomorphisms, inducing the
usual Lusztig's symmetries as $\Bbb Q(q)$-automorphisms not only on
the standard quantum group $U_q(G_2)$ of Drinfel'd-Jimbo type but
also on the centralizd quantum group $U_q^c(G_2)$, only when
$r=q=s^{-1}$. \hfill $\Box$
\end{theorem}

\bigskip
\section{Appendix: Some Calculations in the proof of Proposition 2.3}
\medskip

For $X = f_{1}^{4}f_2$, the relevant terms of $\Delta^{(4)}(X)$ in
$(\ast)$ and their paring-values are as follows:

\hspace*{\fill} TABULAR 1 \hspace*{\fill}

\begin{center}
\begin{tabular}{|l|p{2.5em}|}\hline
\multicolumn{1}{|c|}{\textbf{SUMMANDS}}&\multicolumn{1}{|c|}{\textbf{1}}
\\\hline

$ f_1\omega_1'^3\omega_2' \otimes f_1\omega_1'^2\omega_2' \otimes
f_1\omega_1'\omega_2' \otimes f_1\omega_2' \otimes f_2$ & $a$
\\\hline
$\omega_1'f_1\omega_1'^2\omega_2' \otimes f_1\omega_1'^2\omega_2'
\otimes f_1\omega_1'\omega_2' \otimes f_1\omega_2' \otimes f_2$ &
$xa$
\\\hline
 $\omega_1'^{2}f_1\omega_1'\omega_2'
\otimes f_1\omega_1'^2\omega_2' \otimes f_1\omega_1'\omega_2'
\otimes f_1\omega_2' \otimes f_2$ & $x^{2}a $
\\\hline
$\omega_1'^{3}f_1\omega_2' \otimes f_1\omega_1'^2\omega_2' \otimes
f_1\omega_1'\omega_2' \otimes f_1\omega_2' \otimes f_2$ & $x^{3}a$
\\\hline

$f_1\omega_1'^3\omega_2' \otimes \omega_1'f_1\omega_1'\omega_2'
\otimes f_1\omega_1'\omega_2' \otimes f_1\omega_2' \otimes f_2$ &
$xa$
\\\hline
$\omega_1'f_1\omega_1'^2\omega_2' \otimes
\omega_1'f_1\omega_1'\omega_2' \otimes f_1\omega_1'\omega_2'
\otimes f_1\omega_2' \otimes f_2$ & $x^{2}a$
\\\hline
$\omega_1'^{2}f_1\omega_1'\omega_2' \otimes
\omega_1'f_1\omega_1'\omega_2' \otimes f_1\omega_1'\omega_2'
\otimes f_1\omega_2' \otimes f_2$ & $x^{3}a$
\\\hline
$\omega_1'^{3}f_1\omega_2' \otimes \omega_1'f_1\omega_1'\omega_2'
\otimes f_1\omega_1'\omega_2' \otimes f_1\omega_2' \otimes f_2$ &
$x^{4}a$
\\\hline

$ f_1\omega_1'^3\omega_2' \otimes \omega_1'^{2}f_1\omega_2'
\otimes f_1\omega_1'\omega_2' \otimes f_1\omega_2' \otimes f_2$ &
$x^{2}a$
\\\hline
$\omega_1'f_1\omega_1'^2\omega_2' \otimes
\omega_1'^{2}f_1\omega_2' \otimes f_1\omega_1'\omega_2' \otimes
f_1\omega_2' \otimes f_2$ & $x^{3}a$
\\\hline

$\omega_1'^{2}f_1\omega_1'\omega_2' \otimes
\omega_1'^{2}f_1\omega_2' \otimes f_1\omega_1'\omega_2' \otimes
f_1\omega_2' \otimes f_2$ & $x^{4}a$
\\\hline
$\omega_1'^{3}f_1\omega_2' \otimes \omega_1'^{2}f_1\omega_2'
\otimes f_1\omega_1'\omega_2' \otimes f_1\omega_2' \otimes f_2$ &
$x^{5}a$
\\\hline

$f_1\omega_1'^3\omega_2' \otimes f_1\omega_1'^{2}\omega_2' \otimes
\omega_1'f_1\omega_2' \otimes f_1\omega_2' \otimes f_2$ & $xa$
\\\hline
$\omega_1'f_1\omega_1'^2\omega_2' \otimes
f_1\omega_1'^{2}\omega_2' \otimes \omega_1'f_1\omega_2' \otimes
f_1\omega_2' \otimes f_2$ & $x^{2}a$
\\\hline
$\omega_1'^{2}f_1\omega_1'\omega_2' \otimes
f_1\omega_1'^{2}\omega_2' \otimes \omega_1'f_1\omega_2' \otimes
f_1\omega_2' \otimes f_2$ & $x^{3}a$
\\\hline
$\omega_1'^{3}f_1\omega_2' \otimes f_1\omega_1'^{2}\omega_2'
\otimes \omega_1'f_1\omega_2' \otimes f_1\omega_2' \otimes f_2$ &
$x^{4}a$
\\\hline


$f_1\omega_1'^3\omega_2' \otimes \omega_1'f_1\omega_1'\omega_2'
\otimes \omega_1'f_1\omega_2' \otimes f_1\omega_2' \otimes f_2$ &
$x^{2}a$
\\\hline
$ \omega_1'f_1\omega_1'^2\omega_2' \otimes
\omega_1'f_1\omega_1'\omega_2' \otimes \omega_1'f_1\omega_2'
\otimes f_1\omega_2' \otimes f_2$ & $x^{3}a$
\\\hline
$\omega_1'^{2}f_1\omega_1'\omega_2' \otimes
\omega_1'f_1\omega_1'\omega_2' \otimes \omega_1'f_1\omega_2'
\otimes f_1\omega_2' \otimes f_2$ & $x^{4}a$
\\\hline
$\omega_1'^{3}f_1\omega_2' \otimes \omega_1'f_1\omega_1'\omega_2'
\otimes \omega_1'f_1\omega_2' \otimes f_1\omega_2' \otimes f_2$ &
$x^{5}a$
\\\hline
$f_1\omega_1'^3\omega_2' \otimes \omega_1'^{2}f_1\omega_2' \otimes
\omega_1'f_1\omega_2' \otimes f_1\omega_2' \otimes f_2$ & $x^{3}a$
\\\hline
$\omega_1'f_1\omega_1'^2\omega_2' \otimes
\omega_1'^{2}f_1\omega_2' \otimes \omega_1'f_1\omega_2' \otimes
f_1\omega_2' \otimes f_2$ & $x^{4}a$
\\\hline
$\omega_1'^{2}f_1\omega_1'\omega_2' \otimes
\omega_1'^{2}f_1\omega_2' \otimes \omega_1'f_1\omega_2' \otimes
f_1\omega_2' \otimes f_2$ & $x^{5}a$
\\\hline
$\omega_1'^{3}f_1\omega_2' \otimes \omega_1'^{2}f_1\omega_2'
\otimes \omega_1'f_1\omega_2' \otimes f_1\omega_2' \otimes f_2$ &
$x^{6}a$
\\\hline
\end{tabular}
\end{center}

\begin{center}
 \large
\begin{tabular}{|l|p{2.5em}|}\hline
\multicolumn{1}{|c|}{\textbf{SUMMANDS}}&\multicolumn{1}{|c|}{\textbf{2}}
\\\hline

$f_1\omega_1'^3\omega_2' \otimes f_1\omega_1'^2\omega_2' \otimes
f_1\omega_1'\omega_2' \otimes \omega_1'f_2 \otimes f_1$  & $\bar
xa$
\\\hline
$\omega_1'f_1\omega_1'^2\omega_2' \otimes f_1\omega_1'^2\omega_2'
\otimes f_1\omega_1'\omega_2' \otimes \omega_1'f_2 \otimes f_1$ &
$\bar xxa$
\\\hline

$\omega_1'^{2}f_1\omega_1'\omega_2' \otimes
f_1\omega_1'^2\omega_2' \otimes f_1\omega_1'\omega_2' \otimes
\omega_1'f_2 \otimes f_1$ & $\bar xx^{2}a$
\\\hline
$\omega_1'^{3}f_1\omega_2' \otimes f_1\omega_1'^2\omega_2' \otimes
f_1\omega_1'\omega_2' \otimes \omega_1'f_2 \otimes f_1$ & $\bar
xx^{3}a$
\\\hline

$f_1\omega_1'^3\omega_2'\otimes \omega_1'f_1\omega_1'\omega_2'
\otimes f_1\omega_1'\omega_2' \otimes \omega_1'f_2 \otimes f_1$ &
$\bar xxa$
\\\hline
$\omega_1'f_1\omega_1'^2\omega_2' \otimes
\omega_1'f_1\omega_1'\omega_2' \otimes f_1\omega_1'\omega_2'
\otimes \omega_1'f_2 \otimes f_1$& $\bar xx^{2}a$
\\\hline
$\omega_1'^{2}f_1\omega_1'\omega_2' \otimes
\omega_1'f_1\omega_1'\omega_2' \otimes f_1\omega_1'\omega_2'
\otimes \omega_1'f_2 \otimes f_1$ & $\bar xx^{3}a$
\\\hline
$\omega_1'^{3}f_1\omega_2' \otimes \omega_1'f_1\omega_1'\omega_2'
\otimes f_1\omega_1'\omega_2' \otimes \omega_1'f_2 \otimes f_1$&
$\bar xx^{4}a$
\\\hline

$ f_1\omega_1'^3\omega_2' \otimes \omega_1'^{2}f_1\omega_2'
\otimes f_1\omega_1'\omega_2' \otimes \omega_1'f_2 \otimes f_1$&
$\bar xx^{2}a$
\\\hline
$\omega_1'f_1\omega_1'^2\omega_2' \otimes
\omega_1'^{2}f_1\omega_2' \otimes f_1\omega_1'\omega_2' \otimes
\omega_1'f_2 \otimes f_1$& $\bar xx^{3}a$
\\\hline
$\omega_1'^{2}f_1\omega_1'\omega_2' \otimes
\omega_1'^{2}f_1\omega_2' \otimes f_1\omega_1'\omega_2' \otimes
\omega_1'f_2 \otimes f_1$ & $\bar xx^{4}a$
\\\hline
$\omega_1'^{3}f_1\omega_2' \otimes \omega_1'^{2}f_1\omega_2'
\otimes f_1\omega_1'\omega_2' \otimes \omega_1'f_2 \otimes f_1$&
$\bar xx^{5}a$
\\\hline
\end{tabular}
\end{center}

\begin{center}
 \large
\begin{tabular}{|l|p{2.5em}|}\hline
\multicolumn{1}{|c|}{\textbf{SUMMANDS}}&\multicolumn{1}{|c|}{\textbf{2}}
\\\hline
$f_1\omega_1'^3\omega_2'\otimes f_1\omega_1'^{2}\omega_2' \otimes
\omega_1'f_1\omega_2' \otimes \omega_1'f_2 \otimes f_1$& $\bar
x{x}a$
\\\hline
$\omega_1'f_1\omega_1'^2\omega_2' \otimes
f_1\omega_1'^{2}\omega_2' \otimes \omega_1'f_1\omega_2' \otimes
\omega_1'f_2 \otimes f_1$ & $\bar xx^{2}a$
\\\hline
$\omega_1'^{2}f_1\omega_1'\omega_2' \otimes
f_1\omega_1'^{2}\omega_2' \otimes \omega_1'f_1\omega_2' \otimes
\omega_1'f_2 \otimes f_1$  & $\bar xx^{3}a$
\\\hline
$\omega_1'^{3}f_1\omega_2' \otimes f_1\omega_1'^{2}\omega_2'
\otimes \omega_1'f_1\omega_2' \otimes \omega_1'f_2 \otimes f_1$
&$\bar xx^{4}a$
\\\hline

$f_1\omega_1'^3\omega_2' \otimes\omega_1'f_1\omega_1'\omega_2'
\otimes \omega_1'f_1\omega_2' \otimes \omega_1'f_2 \otimes f_1$ &
$\bar xx^{2}a$
\\\hline

$ \omega_1'f_1\omega_1'^2\omega_2' \otimes
\omega_1'f_1\omega_1'\omega_2' \otimes \omega_1'f_1\omega_2'
\otimes \omega_1'f_2 \otimes f_1$ & $\bar xx^{3}a$
\\\hline
$\omega_1'^{2}f_1\omega_1'\omega_2' \otimes
\omega_1'f_1\omega_1'\omega_2' \otimes \omega_1'f_1\omega_2'
\otimes \omega_1'f_2 \otimes f_1$ &$\bar xx^{4}a$
\\
\hline

\end{tabular}
\end{center}

\begin{center}
 \large
\begin{tabular}{|l|p{2.5em}|}\hline
\multicolumn{1}{|c|}{\textbf{SUMMANDS}}&\multicolumn{1}{|c|}{\textbf{3}}
\\\hline

$f_1\omega_1'^3\omega_2' \otimes f_1\omega_1'^2\omega_2' \otimes
\omega_1'^{2}f_2 \otimes \omega_1'f_1 \otimes f_1$ &$\bar x^{2}xa$
\\\hline
$\omega_1'f_1\omega_1'^2\omega_2' \otimes f_1\omega_1'^2\omega_2'
\otimes \omega_1'^{2}f_2 \otimes \omega_1'f_1 \otimes f_1$ &$\bar
x^{2}x^{2}a$
\\\hline
 $\omega_1'^{2}f_1\omega_1'\omega_2'
\otimes f_1\omega_1'^2\omega_2' \otimes \omega_1'^{2}f_2 \otimes
\omega_1'f_1 \otimes f_1$ & $\bar x^{2}x^{3}a$
\\\hline
 $\omega_1'^{3}f_1\omega_2' \otimes f_1\omega_1'^2\omega_2'
\otimes \omega_1'^{2}f_2 \otimes \omega_1'f_1 \otimes f_1$ &$\bar
x^{2}x^{4}a$
\\\hline
$f_1\omega_1'^3\omega_2' \otimes \omega_1'f_1\omega_1'\omega_2'
\otimes \omega_1'^{2}f_2 \otimes \omega_1'f_1 \otimes f_1$ &$\bar
x^{2}x^{2}a$
\\\hline
$\omega_1'f_1\omega_1'^2\omega_2' \otimes
\omega_1'f_1\omega_1'\omega_2' \otimes \omega_1'^{2}f_2 \otimes
\omega_1'f_1 \otimes f_1$ & $\bar x^{2}x^{3}a$
\\\hline

$\omega_1'^{2}f_1\omega_1'\omega_2' \otimes
\omega_1'f_1\omega_1'\omega_2' \otimes \omega_1'^{2}f_2 \otimes
\omega_1'f_1 \otimes f_1$ & $\bar x^{2}x^{4}a$
\\
\hline

\end{tabular}
\end{center}


\begin{center}
\large
\begin{tabular}{|l|p{2.5em}|}\hline
\multicolumn{1}{|c|}{\textbf{SUMMANDS}}&\multicolumn{1}{|c|}{\textbf{4}}
\\\hline

$f_1\omega_1'^3\omega_2' \otimes \omega_1'^{3}f_2 \otimes
f_1\omega_1'^2 \otimes \omega_1'f_1 \otimes f_1$ & $\bar x^{3}xa$
\\\hline
$\omega_1'f_1\omega_1'^2\omega_2' \otimes \omega_1'^{3}f_2 \otimes
f_1\omega_1'^2 \otimes \omega_1'f_1 \otimes f_1$ & $\bar
x^{3}x^{2}a$
 \\\hline
$\omega_1'^{2}f_1\omega_1'\omega_2' \otimes \omega_1'^{3}f_2
\otimes f_1\omega_1'^2 \otimes \omega_1'f_1 \otimes f_1$ & $\bar
x^{3}x^{3}a$
 \\\hline

$\omega_1'^{3}f_1\omega_2' \otimes \omega_1'^{3}f_2 \otimes
f_1\omega_1'^2 \otimes \omega_1'f_1 \otimes f_1$& $\bar
x^{3}x^{4}a$
\\\hline
\end{tabular}
\end{center}


\begin{center}
\large
\begin{tabular}{|l|p{2.5em}|}\hline
\multicolumn{1}{|c|}{\textbf{SUMMANDS}}&\multicolumn{1}{|c|}{\textbf{4}}
\\\hline
$f_1\omega_1'^3\omega_2' \otimes \omega_1'^{3}f_2 \otimes
\omega_1'f_1\omega_1' \otimes \omega_1'f_1 \otimes f_1$ & $\bar
x^{3}x^{2}a$
\\\hline

$\omega_1'f_1\omega_1'^2\omega_2'\otimes \omega_1'^{3}f_2 \otimes
\omega_1'f_1\omega_1' \otimes \omega_1'f_1 \otimes f_1$& $\bar
x^{3}x^{3}a $
\\\hline

$\omega_1'^{2}f_1\omega_1'\omega_2' \otimes \omega_1'^{3}f_2
\otimes \omega_1'f_1\omega_1' \otimes \omega_1'f_1 \otimes f_1$ &$
\bar x^{3}x^{4}a$
 \\\hline
$\omega_1'^{3}f_1\omega_2' \otimes \omega_1'^{3}f_2 \otimes
\omega_1'f_1\omega_1' \otimes \omega_1'f_1 \otimes f_1$& $\bar
x^{3}x^{5}a$
\\\hline
$f_1\omega_1'^3\omega_2' \otimes \omega_1'^{3}f_2 \otimes
\omega_1'^{2}f_1 \otimes \omega_1'f_1 \otimes f_1$& $\bar
x^{3}x^{3}a$
\\\hline
$\omega_1'f_1\omega_1'^2\omega_2'\otimes \omega_1'^{3}f_2 \otimes
\omega_1'^{2}f_1 \otimes \omega_1'f_1 \otimes f_1$& $\bar
x^{3}x^{4}a $
\\\hline
$\omega_1'^{2}f_1\omega_1'\omega_2' \otimes \omega_1'^{3}f_2
\otimes \omega_1'^{2}f_1\omega_1' \otimes \omega_1'f_1 \otimes
f_1$ &$ \bar x^{3}x^{5}a$
 \\
\hline

\end{tabular}
\end{center}

\begin{center}
\large
\begin{tabular}{|l|p{2.5em}|}\hline
\multicolumn{1}{|c|}{\textbf{SUMMANDS}}&\multicolumn{1}{|c|}{\textbf{5}}
\\\hline

$\omega_1'^{4}f_2 \otimes f_1\omega_1'^3 \otimes f_1\omega_1'^2
\otimes \omega_1'f_1 \otimes f_1$ & $\bar x^{4}xa$
\\\hline
$\omega_1'^{4}f_2 \otimes \omega_1'f_1\omega_1'^2 \otimes
f_1\omega_1'^2 \otimes \omega_1'f_1 \otimes f_1$  & $\bar
x^{4}x^{2}a$
\\\hline
$\omega_1'^{4}f_2 \otimes \omega_1'^{2}f_1\omega_1' \otimes
f_1\omega_1'^2 \otimes \omega_1'f_1 \otimes f_1$  & $\bar
x^{4}x^{3}a$
\\\hline
$\omega_1'^{4}f_2 \otimes \omega_1'^{3}f_1 \otimes f_1\omega_1'^2
\otimes \omega_1'f_1 \otimes f_1$& $\bar x^{4}x^{4}a$
\\\hline

$\omega_1'^{4}f_2 \otimes f_1\omega_1'^3 \otimes
\omega_1'f_1\omega_1' \otimes \omega_1'f_1 \otimes f_1$ &$\bar
x^{4}x^{2}a$
\\\hline
$\omega_1'^{4}f_2 \otimes \omega_1'f_1\omega_1'^2 \otimes
\omega_1'f_1\omega_1' \otimes \omega_1'f_1 \otimes f_1$& $\bar
x^{4}x^{3}a$
\\\hline
$\omega_1'^{4}f_2 \otimes \omega_1'^{2}f_1\omega_1' \otimes
\omega_1'f_1\omega_1' \otimes \omega_1'f_1 \otimes f_1$& $\bar
x^{4}x^{4}a$
\\
\hline
\end{tabular}
\end{center}

\begin{center}
\large
\begin{tabular}{|l|p{2.5em}|}\hline
\multicolumn{1}{|c|}{\textbf{SUMMANDS}}&\multicolumn{1}{|c|}{\textbf{5}}
\\\hline
$\omega_1'^{4}f_2 \otimes f_1\omega_1'^3 \otimes \omega_1'^{2}f_1
\otimes f_1\omega_1' \otimes f_1$&$\bar x^{4}x^{2}a$
\\\hline
$\omega_1'^{4}f_2 \otimes \omega_1'f_1\omega_1'^2 \otimes
\omega_1'^{2}f_1 \otimes 'f_1\omega_1 \otimes f_1$ & $\bar
x^{4}x^{3}a$
\\\hline

$\omega_1'^{4}f_2 \otimes \omega_1'^{2}f_1\omega_1' \otimes
\omega_1'^{2}f_1 \otimes f_1\omega_1' \otimes f_1$ & $\bar
x^{4}x^{4}a$
\\\hline
$\omega_1'^{4}f_2 \otimes \omega_1'^{3}f_1 \otimes
\omega_1'^{2}f_1 \otimes f_1\omega_1' \otimes f_1$ &$\bar
x^{4}x^{5}a$
\\\hline
\end{tabular}
\end{center}

The relevant terms and the pairing-values of $\Delta^{(4)}(X)$ for
$X = f_{2}f_{1}^{4}$ in $(\ast)$ are as follows:

\hspace*{\fill} TABULAR 2 \hspace*{\fill}

\begin{center}
\large
\begin{tabular}{|l|p{2.5em}|}\hline
\multicolumn{1}{|c|}{\textbf{SUMMANDS}}&\multicolumn{1}{|c|}{\textbf{1}}
\\\hline

$ \omega_2'f_1\omega_1'^3 \otimes \omega_2'f_1\omega_1'^2 \otimes
\omega_2'f_1\omega_1' \otimes \omega_2'f_1\otimes f_2$ & $ay^4$
\\\hline
$\omega_1'\omega_2'f_1\omega_1'^2 \otimes \omega_2'f_1\omega_1'^2
\otimes \omega_2'f_1\omega_1' \otimes \omega_2'f_1\otimes f_2$ &
$axy^4$
\\\hline
 $\omega_1'^{2}\omega_2'f_1\omega_1'
\otimes \omega_2'f_1\omega_1'^2 \otimes \omega_2'f_1\omega_1'
\otimes \omega_2'f_1\otimes f_2$ & $ax^{2}y^4$
\\\hline
$\omega_1'^{3}\omega_2'f_1 \otimes \omega_2'f_1\omega_1'^2 \otimes
\omega_2'f_1\omega_1' \otimes \omega_2'f_1\otimes f_2$ &
$ax^{3}y^4$
\\\hline

$\omega_2'f_1\omega_1'^3 \otimes \omega_1'\omega_2'f_1\omega_1'
\otimes \omega_2'f_1\omega_1' \otimes \omega_2'f_1\otimes f_2$ &
$axy^4$
\\\hline
$\omega_1'\omega_2'f_1\omega_1'^2 \otimes
\omega_1'\omega_2'f_1\omega_1' \otimes \omega_2'f_1\omega_1'
\otimes \omega_2'f_1\otimes f_2$ & $ax^{2}y^4$
\\\hline
$\omega_1'^{2}\omega_2'f_1\omega_1'\otimes
\omega_1'\omega_2'f_1\omega_1'\otimes \omega_2'f_1\omega_1'
\otimes \omega_2'f_1\otimes f_2$ & $ax^{3}y^4$
\\
\end{tabular}
\end{center}

\begin{center}
\large
\begin{tabular}{|l|p{2.5em}|}\hline
\multicolumn{1}{|c|}{\textbf{SUMMANDS}}&\multicolumn{1}{|c|}{\textbf{2}}
\\\hline
$\omega_2'f_1\omega_1'^3 \otimes \omega_2'f_1\omega_1'^2 \otimes
\omega_2'f_1\omega_1' \otimes f_2\omega_1' \otimes f_1$  & $ay^3$
\\\hline
$\omega_1'\omega_2'f_1\omega_1'^2 \otimes \omega_2'f_1\omega_1'^2
\otimes \omega_2'f_1\omega_1' \otimes f_2\omega_1' \otimes f_1$ &
 $axy^3$
\\\hline

$\omega_1'^{2}\omega_2'f_1\omega_1' \otimes
\omega_2'f_1\omega_1'^2 \otimes \omega_2'f_1\omega_1' \otimes
f_2\omega_1' \otimes f_1$ & $ax^{2}y^3$
\\\hline

$\omega_1'^{3}\omega_2'f_1 \otimes \omega_2'f_1\omega_1'^2 \otimes
\omega_2'f_1\omega_1' \otimes f_2\omega_1' \otimes f_1$ &
$ax^{3}y^3$
\\\hline

$\omega_2'f_1\omega_1'^3\otimes \omega_1'\omega_2'f_1\omega_1'
\otimes \omega_2'f_1\omega_1' \otimes f_2\omega_1' \otimes f_1$ &
 $axy^3$
\\\hline
$\omega_1'\omega_2'f_1\omega_1'^2 \otimes
\omega_1'\omega_2'f_1\omega_1' \otimes \omega_2'f_1\omega_1'
\otimes f_2\omega_1' \otimes f_1$  & $ax^{2}y^3$
\\\hline
$\omega_1'^{2}\omega_2'f_1\omega_1' \otimes
\omega_1'\omega_2'f_1\omega_1' \otimes \omega_2'f_1\omega_1'
\otimes f_2\omega_1' \otimes f_1$  & $ax^{3}y^3$
\\\hline
$\omega_1'^{3}\omega_2'f_1 \otimes \omega_1'\omega_2'f_1\omega_1'
\otimes \omega_2'f_1\omega_1' \otimes f_2\omega_1' \otimes f_1$ &
 $ax^{4}y^3$
\\\hline
\end{tabular}
\end{center}

\begin{center}
\large
\begin{tabular}{|l|p{2.5em}|}\hline
\multicolumn{1}{|c|}{\textbf{SUMMANDS}}&\multicolumn{1}{|c|}{\textbf{2}}
\\\hline
$ \omega_2'f_1\omega_1'^3 \otimes \omega_1'^{2}\omega_2'f_1
\otimes \omega_2'f_1\omega_1' \otimes f_2\omega_1' \otimes f_1$ &
 $ax^{2}y^3$
\\\hline
$\omega_1'\omega_2'f_1\omega_1'^2 \otimes
\omega_1'^{2}\omega_2'f_1 \otimes \omega_2'f_1\omega_1' \otimes
f_2\omega_1' \otimes f_1$ &  $ax^{3}y^3$
\\\hline
$\omega_1'^{2}\omega_2'f_1\omega_1' \otimes
\omega_1'^{2}\omega_2'f_1 \otimes \omega_2'f_1\omega_1' \otimes
f_2\omega_1' \otimes f_1$ &$ax^{4}y^3$
\\\hline
$\omega_1'^{3}\omega_2'f_1 \otimes \omega_1'^{2}\omega_2'f_1
\otimes \omega_2'f_1\omega_1' \otimes f_2\omega_1' \otimes f_1$ &
 $ax^{5}y^3$
\\\hline

$\omega_2'f_1\omega_1'^3\otimes \omega_2'f_1\omega_1'^{2} \otimes
\omega_1'\omega_2'f_1 \otimes f_2\omega_1' \otimes f_1$ & $axy^3$
\\\hline
$\omega_1'\omega_2'f_1\omega_1'^2 \otimes
\omega_2'f_1\omega_1'^{2} \otimes \omega_1'\omega_2'f_1 \otimes
f_2\omega_1' \otimes f_1$ & $ax^{2}y^3$
\\\hline
$\omega_1'^{2}\omega_2'f_1\omega_1' \otimes
\omega_2'f_1\omega_1'^{2} \otimes \omega_1'\omega_2'f_1 \otimes
f_2\omega_1' \otimes f_1$ & $ax^{3}y^3$
\\
\hline

\end{tabular}
\end{center}

\begin{center}
\large
\begin{tabular}{|l|p{2.5em}|}\hline
\multicolumn{1}{|c|}{\textbf{SUMMANDS}}&\multicolumn{1}{|c|}{\textbf{3}}
\\\hline

$\omega_2'f_1\omega_1'^3\otimes \omega_2'f_1\omega_1'^2\otimes
f_2\omega_1'^{2} \otimes \omega_1'f_1 \otimes f_1$ &$axy^{2}$
\\\hline
$\omega_1'\omega_2'f_1\omega_1'^2
\otimes\omega_2'f_1\omega_1'^2\otimes f_2\omega_1'^{2} \otimes
\omega_1'f_1 \otimes f_1$ & $ax^{2}y^{2}$
\\\hline

$\omega_1'^{2}\omega_2'f_1\omega_1'\otimes
\omega_2'f_1\omega_1'^2\otimes f_2\omega_1'^{2} \otimes
\omega_1'f_1 \otimes f_1$ &$ax^{3}y^{2}$
\\\hline

$\omega_1'^{3}\omega_2'f_1 \otimes \omega_2'f_1\omega_1'^2\otimes
f_2\omega_1'^{2} \otimes \omega_1'f_1 \otimes f_1$ &$ax^{4}y^{2}$
\\\hline

$\omega_2'f_1\omega_1'^3 \otimes \omega_1'\omega_2'f_1\omega_1'
\otimes f_2\omega_1'^{2} \otimes \omega_1'f_1 \otimes f_1$
&$ax^{2}y^{2}$
\\\hline
$\omega_1'\omega_2'f_1\omega_1'^2 \otimes
\omega_1'\omega_2'f_1\omega_1' \otimes f_2\omega_1'^{2} \otimes
\omega_1'f_1 \otimes f_1$  &$ax^{3}y^{2}$
\\\hline

$\omega_1'^{2}\omega_2'f_1\omega_1' \otimes
\omega_1'\omega_2'f_1\omega_1' \otimes f_2\omega_1'^{2} \otimes
\omega_1'f_1 \otimes f_1$ & $ax^{4}y^{2}$
\\
\hline
\end{tabular}
\end{center}

\begin{center}
\large
\begin{tabular}{|l|p{2.5em}|}\hline
\multicolumn{1}{|c|}{\textbf{SUMMANDS}}&\multicolumn{1}{|c|}{\textbf{4}}
\\\hline

$\omega_2'f_1\omega_1'^3 \otimes f_2\omega_1'^{3} \otimes
f_1\omega_1'^2 \otimes \omega_1'f_1 \otimes f_1$ & $axy$
\\\hline
$\omega_1'\omega_2'f_1\omega_1'^2 \otimes f_2\omega_1'^{3} \otimes
f_1\omega_1'^2 \otimes \omega_1'f_1 \otimes f_1$ & $ax^{2}y$
\\\hline
$\omega_1'^{2}\omega_2'f_1\omega_1' \otimes f_2\omega_1'^{3}
\otimes f_1\omega_1'^2 \otimes \omega_1'f_1 \otimes f_1$
&$ax^{3}y$
\\\hline
$\omega_1'^{3}\omega_2'f_1 \otimes f_2\omega_1'^{3} \otimes
f_1\omega_1'^2 \otimes \omega_1'f_1 \otimes f_1$  & $ax^{4}y$
\\\hline
$\omega_2'f_1\omega_1'^3 \otimes f_2\omega_1'^{3} \otimes
\omega_1'f_1\omega_1' \otimes \omega_1'f_1 \otimes f_1$ &
$ax^{2}y$
\\\hline
$\omega_1'\omega_2'f_1\omega_1'^2\otimes f_2\omega_1'^{3} \otimes
\omega_1'f_1\omega_1' \otimes \omega_1'f_1 \otimes f_1$& $ax^{3}y$
\\\hline

$\omega_1'^{2}\omega_2'f_1\omega_1' \otimes f_2\omega_1'^{3}
\otimes \omega_1'f_1\omega_1' \otimes \omega_1'f_1 \otimes f_1$&
$ax^{4}y$
\\
\hline
\end{tabular}
\end{center}

\begin{center}
\large
\begin{tabular}{|l|p{2.5em}|}\hline
\multicolumn{1}{|c|}{\textbf{SUMMANDS}}&\multicolumn{1}{|c|}{\textbf{5}}
\\\hline
$f_2\omega_1'^{4} \otimes f_1\omega_1'^3 \otimes f_1\omega_1'^2
\otimes \omega_1'f_1 \otimes f_1$ & $ax$
\\\hline
$f_2\omega_1'^{4} \otimes \omega_1'f_1\omega_1'^2 \otimes
f_1\omega_1'^2 \otimes \omega_1'f_1 \otimes f_1$ & $ax^2$
\\\hline
$f_2\omega_1'^{4} \otimes \omega_1'^{2}f_1\omega_1' \otimes
f_1\omega_1'^2 \otimes \omega_1'f_1 \otimes f_1$ & $ax^{3}$
\\\hline
$f_2\omega_1'^{4} \otimes \omega_1'^{3}f_1 \otimes f_1\omega_1'^2
\otimes \omega_1'f_1 \otimes f_1$ &  $ax^{4}$
\\\hline

$f_2\omega_1'^{4} \otimes f_1\omega_1'^3 \otimes
\omega_1'f_1\omega_1' \otimes \omega_1'f_1 \otimes f_1$ &$ax^{2}$
\\\hline
$f_2\omega_1'^{4} \otimes \omega_1'f_1\omega_1'^2 \otimes
\omega_1'f_1\omega_1' \otimes \omega_1'f_1 \otimes f_1$ & $ax^{3}$
\\\hline
$f_2\omega_1'^{4} \otimes \omega_1'^{2}f_1\omega_1' \otimes
\omega_1'f_1\omega_1' \otimes \omega_1'f_1 \otimes f_1$ & $ax^{4}$
\\
\hline
\end{tabular}
\end{center}

\begin{center}
\large
\begin{tabular}{|l|p{2.5em}|}\hline
\multicolumn{1}{|c|}{\textbf{SUMMANDS}}&\multicolumn{1}{|c|}{\textbf{5}}
\\\hline
$f_2\omega_1'^{4} \otimes f_1\omega_1'^3 \otimes
\omega_1'f_1\omega_1' \otimes f_1\omega_1' \otimes f_1$ &$ax$
\\\hline
$f_2\omega_1'^{4} \otimes \omega_1'f_1\omega_1'^2 \otimes
\omega_1'f_1\omega_1' \otimes f_1\omega_1' \otimes f_1$ & $ax^{2}$
\\\hline
$f_2\omega_1'^{4} \otimes \omega_1'^{2}f_1\omega_1' \otimes
\omega_1'f_1\omega_1' \otimes f_1\omega_1' \otimes f_1$& $ax^{3}$
\\\hline
$f_2\omega_1'^{4} \otimes \omega_1'^{3}f_1 \otimes
\omega_1'f_1\omega_1' \otimes f_1\omega_1' \otimes f_1$&$ax^{4}$
\\\hline
$f_2\omega_1'^{4} \otimes f_1\omega_1'^3 \otimes \omega_1'^{2}f_1
\otimes f_1\omega_1' \otimes f_1$ &$ax^{2}$
\\\hline
$f_2\omega_1'^{4} \otimes \omega_1'f_1\omega_1'^2 \otimes
\omega_1'^{2}f_1 \otimes f_1\omega_1' \otimes f_1$ & $ax^{3}$
\\\hline
$f_2\omega_1'^{4} \otimes \omega_1'^{2}f_1\omega_1' \otimes
\omega_1'^{2}f_1 \otimes f_1\omega_1' \otimes f_1$ & $ax^{4}$
\\\hline
$f_2\omega_1'^{4} \otimes \omega_1'^{3}f_1 \otimes
\omega_1'^{2}f_1 \otimes f_1\omega_1' \otimes f_1$ &$ax^{5}$
\\\hline
\end{tabular}
\end{center}

The relevant terms and the pairing-values of $\Delta^{(4)}(X)$ for
$X = f_{1}^{2}f_{2}f_{1}^{2}$ in $(\ast)$ are as follows:

\hspace*{\fill} TABULAR 3 \hspace*{\fill}
\begin{center}
 \large
\begin{tabular}{|l|p{2.5em}|}\hline
\multicolumn{1}{|c|}{\textbf{SUMMANDS}}&\multicolumn{1}{|c|}{\textbf{1}}
\\\hline

$ f_1\omega_1'^3\omega_2' \otimes f_1\omega_1'^2\omega_2' \otimes
\omega_2'f_1\omega_1' \otimes \omega_2'f_1\otimes f_2$ & $ay^2$
\\\hline
$\omega_1'f_1\omega_1'^2\omega_2' \otimes f_1\omega_1'^2\omega_2'
\otimes \omega_2'f_1\omega_1' \otimes \omega_2'f_1\otimes f_2$ &
$axy^2$
\\\hline
$f_1\omega_1'^{3}\omega_2' \otimes f_1\omega_1'^2\omega_2' \otimes
\omega_1'\omega_2'f_1 \otimes \omega_2'f_1\otimes f_2$ & $axy^2$
\\\hline
$\omega_1'f_1\omega_1'^{2}\omega_2' \otimes
f_1\omega_1'^2\omega_2' \otimes \omega_1'\omega_2'f_1 \otimes
\omega_2'f_1\otimes f_2$ & $ax^{2}y^2$
\\\hline

$f_1\omega_1'^3\omega_2' \otimes \omega_1'\omega_2'f_1\omega_1'
\otimes f_1\omega_1'\omega_2' \otimes \omega_2'f_1\otimes f_2$ &
$axy^2$
\\\hline
$\omega_1'f_1\omega_1'^2\omega_2' \otimes
\omega_1'\omega_2'f_1\omega_1' \otimes f_1\omega_1'\omega_2'
\otimes \omega_2'f_1\otimes f_2$ & $ax^{2}y^2$
\\\hline
$f_1\omega_1'^{3}\omega_2'\otimes \omega_1'^2\omega_2'f_1\otimes
f_1\omega_1'\omega_2' \otimes \omega_2'f_1\otimes f_2$ &
$ax^{2}y^2$
\\
\hline
\end{tabular}
\end{center}

\begin{center}
 \large
\begin{tabular}{|l|p{2.5em}|}\hline
\multicolumn{1}{|c|}{\textbf{SUMMANDS}}&\multicolumn{1}{|c|}{\textbf{2}}
\\\hline
$f_1\omega_1'^{3}\omega_2' \otimes f_1\omega_1'^{2}\omega_2'
 \otimes \omega_2'f_1\omega_1' \otimes f_2\omega_1'\otimes
f_1$ & $ay$
\\\hline
$f_1\omega_1'^{3}\omega_2' \otimes f_1\omega_1'^{2}\omega_2'
 \otimes \omega_1'\omega_2'f_1 \otimes f_2\omega_1'\otimes
f_1$& $axy$
\\\hline
$f_1\omega_1'^{3}\omega_2' \otimes \omega_1'\omega_2'f_1\omega_1'
 \otimes f_1\omega_1'\omega_2' \otimes f_2\omega_1'\otimes
f_1$ & $axy$
\\\hline
\end{tabular}
\end{center}

\begin{center}
 \large
\begin{tabular}{|l|p{2.5em}|}\hline
\multicolumn{1}{|c|}{\textbf{SUMMANDS}}&\multicolumn{1}{|c|}{\textbf{2}}
\\\hline
$f_1\omega_1'^{3}\omega_2' \otimes \omega_1'^2\omega_2'f_1
 \otimes f_1\omega_1'\omega_2' \otimes f_2\omega_1'\otimes
f_1$ & $ax^{2}y$
\\\hline

$\omega_1'f_1\omega_1'^{2}\omega_2' \otimes
f_1\omega_1'^2\omega_2'
 \otimes \omega_2'f_1\omega_1' \otimes f_2\omega_1'\otimes
f_1$& $axy$
\\\hline
$\omega_1'f_1\omega_1'^{2}\omega_2' \otimes
f_1\omega_1'^2\omega_2'
 \otimes \omega_1'\omega_2'f_1 \otimes f_2\omega_1'\otimes
f_1$ & $ax^{2}y$
\\\hline

$\omega_1'^{2}\omega_2'f_1\omega_1' \otimes
f_1\omega_1'^2\omega_2'
 \otimes f_1\omega_1'\omega_2' \otimes f_2\omega_1'\otimes
f_1$ & $ax^{2}y$
\\\hline
$\omega_1'^{3}\omega_2'f_1 \otimes f_1\omega_1'^2\omega_2'
 \otimes f_1\omega_1'\omega_2' \otimes f_2\omega_1'\otimes
f_1$ &$ax^{3}y$
\\\hline

$\omega_1'f_1\omega_1'^{2}\omega_2' \otimes
\omega_1'\omega_2'f_1\omega_1'
 \otimes f_1\omega_1'\omega_2' \otimes f_2\omega_1'\otimes
f_1$ &$ax^{2}y$
\\\hline
$\omega_1'f_1\omega_1'^{2}\omega_2' \otimes
\omega_1'^2\omega_2'f_1
 \otimes f_1\omega_1'\omega_2' \otimes f_2\omega_1'\otimes
f_1$ & $ax^{3}y$
\\
\hline
\end{tabular}
\end{center}

\begin{center}
 \large
\begin{tabular}{|l|p{2.5em}|}\hline
\multicolumn{1}{|c|}{\textbf{SUMMANDS}}&\multicolumn{1}{|c|}{\textbf{3}}
\\\hline
$f_1\omega_1'^3\omega_2' \otimes f_1\omega_1'^2\omega_2' \otimes
f_2\omega_1'^2 \otimes f_1\omega_1' \otimes f_1$ & a
\\\hline
$ f_1\omega_1'^3\omega_2' \otimes \omega_1'\omega_2'f_1\omega_1'
\otimes \omega_1'f_2\omega_1' \otimes f_1\omega_1'\otimes
f_1$&$ax\bar{x}y$
\\\hline

$f_1\omega_1'^3\omega_2' \otimes \omega_1'^2\omega_2'f_1 \otimes
\omega_1'f_2\omega_1' \otimes f_1\omega_1'\otimes f_1$ &
$ax^2\bar{x}y$
\\\hline
 $\omega_1'f_1\omega_1'^{2}\omega_2'
\otimes f_1\omega_1'^2\omega_2' \otimes f_2\omega_1'^2 \otimes
f_1\omega_1'\otimes f_1$&$ax$
\\\hline
$ \omega_1'^2\omega_2'f_1\omega_1' \otimes f_1\omega_1'^2\omega_2'
\otimes \omega_1'f_2\omega_1' \otimes f_1\omega_1'\otimes f_1$
&$ax^2\bar{x}y$
\\\hline
$ \omega_1'^3\omega_2'f_1 \otimes f_1\omega_1'^2\omega_2' \otimes
\omega_1'f_2\omega_1' \otimes f_1\omega_1'\otimes f_1$
&$ax^3\bar{x}y$
\\\hline

$\omega_1'f_1\omega_1'^2\omega_2' \otimes
\omega_1'\omega_2'f_1\omega_1' \otimes \omega_1'f_2\omega_1'
\otimes f_1\omega_1'\otimes f_1$  & $ax^2\bar{x}y$
\\
\hline
\end{tabular}
\end{center}

\begin{center}
 \large
\begin{tabular}{|l|p{2.5em}|}\hline
\multicolumn{1}{|c|}{\textbf{SUMMANDS}}&\multicolumn{1}{|c|}{\textbf{3}}
\\\hline
$\omega_1'f_1\omega_1'^{2}\omega_2' \otimes
\omega_1'^2\omega_2'f_1 \otimes \omega_1'f_2\omega_1' \otimes
\omega_1'f_1 \otimes f_1$ & $ax^4\bar{x}y$
\\\hline

$\omega_1'^{2}\omega_2'f_1\omega_1' \otimes
\omega_1'f_1\omega_1'\omega_2' \otimes \omega_1'f_2\omega_1'
\otimes \omega_1'f_1 \otimes f_1$  &$ax^4\bar{x}y$
\\\hline

$\omega_1'^{3}\omega_2'f_1 \otimes \omega_1'f_1\omega_1'\omega_2'
\otimes \omega_1'f_2\omega_1' \otimes \omega_1'f_1 \otimes f_1$ &
$ax^5\bar{x}y$
\\\hline
$\omega_1'^{2}\omega_2'f_1\omega_1' \otimes
\omega_1'^{2}\omega_2'f_1
 \otimes \omega_1'^{2}f_2 \otimes \omega_1'f_1 \otimes
f_1$ & $ax^5\bar{x}^{2}y^2$
\\\hline
$\omega_1'^{3}\omega_2'f_1 \otimes \omega_1'^{2}\omega_2'f_1
 \otimes \omega_1'^{2}f_2 \otimes \omega_1'f_1 \otimes
f_1$ & $ax^6\bar{x}^{2}y^2$
\\\hline
\end{tabular}
\end{center}

\begin{center}
 \large
\begin{tabular}{|l|p{2.5em}|}\hline
\multicolumn{1}{|c|}{\textbf{SUMMANDS}}&\multicolumn{1}{|c|}{\textbf{4}}
\\\hline
$f_1\omega_1'^{3}\omega_2' \otimes \omega_1'f_2\omega_1'^2 \otimes
f_1\omega_1'^2 \otimes f_1\omega_1' \otimes f_1$ & $a\bar{x}$
\\\hline

$f_1\omega_1'^{3}\omega_2' \otimes \omega_1'f_2\omega_1'^{2}
\otimes\omega_1'f_1\omega_1' \otimes f_1\omega_1' \otimes f_1$ &
$ax\bar{x}$
\\\hline
$f_1\omega_1'^{3}\omega_2' \otimes \omega_1'f_2\omega_1'^{2}
\otimes \omega_1'^{2}f_1 \otimes f_1\omega_1' \otimes
f_1$&$ax^2\bar{x}$
\\\hline
$ \omega_1'f_1\omega_1'^2\omega_2' \otimes \omega_1'f_2\omega_1'^2
\otimes f_1\omega_1'^2 \otimes f_1\omega_1'\otimes f_1$ &
$ax\bar{x}$
\\\hline
$\omega_1'^2\omega_2'f_1\omega_1' \otimes
\omega_1'^{2}f_2\omega_1' \otimes f_1\omega_1'^{2} \otimes
f_1\omega_1'\otimes f_1$ & $ax^2\bar{x}^{2}y$
\\\hline
$\omega_1'^3\omega_2'f_1 \otimes \omega_1'^{2}f_2\omega_1' \otimes
f_1\omega_1'^{2} \otimes f_1\omega_1'\otimes f_1$ &
$ax^3\bar{x}^{2}y$
\\\hline
$\omega_1'f_1\omega_1'^2\omega_2' \otimes \omega_1'f_2\omega_1'^2
\otimes \omega_1'f_1\omega_1' \otimes f_1\omega_1'\otimes f_1$ &
$ax^2\bar{x}$
\\
\hline
\end{tabular}
\end{center}

\begin{center}
 \large
\begin{tabular}{|l|p{2.5em}|}\hline
\multicolumn{1}{|c|}{\textbf{SUMMANDS}}&\multicolumn{1}{|c|}{\textbf{5}}
\\\hline
$\omega_1'^2f_2\omega_1'^{2} \otimes f_1\omega_1'^3 \otimes
f_1\omega_1'^2 \otimes f_1\omega_1'\otimes f_1$ &$a\bar{x}^{2}$
\\\hline

$\omega_1'^2f_2\omega_1'^{2} \otimes f_1\omega_1'^3 \otimes
\omega_1'f_1\omega_1' \otimes f_1\omega_1'\otimes
f_1$&$ax\bar{x}^{2}$
\\\hline
$\omega_1'^2f_2\omega_1'^{2} \otimes f_1\omega_1'^3 \otimes
\omega_1'^{2}f_1\omega_1' \otimes f_1\omega_1'\otimes f_1$
&$ax^2\bar{x}^{2}$
\\\hline
$\omega_1'^2f_2\omega_1'^{2} \otimes \omega_1'f_1\omega_1'^2
\otimes f_1\omega_1'^{2} \otimes f_1\omega_1'\otimes
f_1$&$ax\bar{x}^{2}$
\\\hline
$\omega_1'^2f_2\omega_1'^{2} \otimes \omega_1'^{2}f_1\omega_1'
\otimes f_1\omega_1'^{2} \otimes f_1\omega_1'\otimes f_1$
&$ax^2\bar{x}^{2}$
\\\hline
$\omega_1'^2f_2\omega_1'^{2} \otimes \omega_1'^{3}f_1 \otimes
f_1\omega_1'^{2} \otimes f_1\omega_1'\otimes f_1$
&$ax^3\bar{x}^{2}$
\\\hline
$\omega_1'^{2}f_2\omega_1'^2\omega_2' \otimes
\omega_1'f_1\omega_1'^2 \otimes \omega_1'f_1\omega_1' \otimes
f_1\omega_1'\otimes f_1$&$ax^2\bar{x}^2$
\\\hline

$\omega_1'^{2}f_2\omega_1'^2\omega_2' \otimes
\omega_1'f_1\omega_1'^2 \otimes \omega_1'^{2}f_1 \otimes
f_1\omega_1'\otimes f_1$ &$ax^3\bar{x}^2$
\\\hline
$\omega_1'^{2}f_2\omega_1'^2 \otimes \omega_1'^{2}f_1\omega_1'
\otimes \omega_1'f_1\omega_1' \otimes f_1\omega_1'\otimes
f_1$&$ax^3\bar{x}^2$
\\\hline
$\omega_1'^{2}f_2\omega_1'^2 \otimes \omega_1'^{3}f_1 \otimes
\omega_1'f_1\omega_1' \otimes f_1\omega_1'\otimes
f_1$&$ax^4\bar{x}^2$
\\\hline
\end{tabular}
\end{center}

\begin{center}
 \large
\begin{tabular}{|l|p{2.5em}|}\hline
\multicolumn{1}{|c|}{\textbf{SUMMANDS}}&\multicolumn{1}{|c|}{\textbf{5}}
\\\hline
$\omega_1'^{2}f_2\omega_1'^2 \otimes \omega_1'^{2}f_1\omega_1'
\otimes \omega_1'^{2}f_1 \otimes f_1\omega_1'\otimes f_1$
&$ax^4\bar{x}^{2}$
\\\hline

$\omega_1'^{2}f_2\omega_1'^2 \otimes \omega_1'^{3}f_1 \otimes
\omega_1'^{2}f_1 \otimes f_1\omega_1'\otimes f_1$
&$ax^5\bar{x}^{2}$
\\\hline
$\omega_1'^{2}f_2\omega_1'^{2} \otimes f_1\omega_1'^3 \otimes
f_1\omega_1'^2 \otimes \omega_1'f_1 \otimes f_1$ &$ax\bar{x}^2$
\\\hline

$\omega_1'^{2}f_2\omega_1'^{2} \otimes f_1\omega_1'^3 \otimes
\omega_1'f_1\omega_1' \otimes \omega_1'f_1 \otimes
f_1$&$ax^2\bar{x}^2$
\\\hline
$\omega_1'^{2}f_2\omega_1'^{2} \otimes f_1\omega_1'^3 \otimes
\omega_1'^{2}f_1 \otimes \omega_1'f_1 \otimes f_1$
&$ax^3\bar{x}^2$
\\\hline
$\omega_1'^{2}f_2\omega_1'^{2} \otimes \omega_1'f_1\omega_1'^2
\otimes f_1\omega_1'^2 \otimes \omega_1'f_1 \otimes
f_1$&$ax^2\bar{x}^2$
\\\hline
$\omega_1'^{2}f_2\omega_1'^{2} \otimes \omega_1'^{2}f_1\omega_1'
\otimes f_1\omega_1'^2 \otimes \omega_1'f_1 \otimes f_1$
&$ax^3\bar{x}^2$
\\
\hline
\end{tabular}
\end{center}

The relevant terms and the pairing-values of $\Delta^{(4)}(X)$ for
$X = f_{1}^{3}f_{2}f_{1}$ in $(\ast)$ are as follows:

\hspace*{\fill} TABULAR 4 \hspace*{\fill}
\begin{center}
 \large
\begin{tabular}{|l|p{2.5em}|}\hline
\multicolumn{1}{|c|}{\textbf{SUMMANDS}}&\multicolumn{1}{|c|}{\textbf{1}}
\\\hline

$ f_1\omega_1'^3\omega_2' \otimes f_1\omega_1'^2\omega_2' \otimes
f_1\omega_2'\omega_1' \otimes \omega_2'f_1\otimes f_2$ & $ay$
\\\hline
$ f_1\omega_1'^3\omega_2' \otimes \omega_1'f_1\omega_1'\omega_2'
\otimes f_1\omega_2'\omega_1' \otimes \omega_2'f_1\otimes f_2$ &
$axy$
\\\hline
$ \omega_1'f_1\omega_1'^2\omega_2' \otimes f_1\omega_1'^2\omega_2'
\otimes f_1\omega_2'\omega_1' \otimes \omega_2'f_1\otimes f_2$ &
$axy$
\\\hline
$ \omega_1'^{2}f_1\omega_1'\omega_2' \otimes
f_1\omega_1'^2\omega_2' \otimes f_1\omega_2'\omega_1' \otimes
\omega_2'f_1\otimes f_2$ & $ax^{2}y$
\\\hline

$ \omega_1'f_1\omega_1'^2\omega_2' \otimes
\omega_1'f_1\omega_1'\omega_2' \otimes f_1\omega_2'\omega_1'
\otimes \omega_2'f_1\otimes f_2$ & $ax^{2}y$
\\\hline
$ \omega_1'^{2}f_1\omega_1'\omega_2' \otimes
\omega_1'f_1\omega_1'\omega_2' \otimes f_1\omega_2'\omega_1'
\otimes \omega_2'f_1\otimes f_2$ & $ax^{3}y$
\\\hline

$ f_1\omega_1'^{3}\omega_2' \otimes f_1\omega_1'^2\omega_2'
\otimes \omega_1'\omega_2'f_1 \otimes f_1\omega_2'\otimes f_2$ &
$axy$
\\
\hline
\end{tabular}
\end{center}

\begin{center}
 \large
\begin{tabular}{|l|p{2.5em}|}\hline
\multicolumn{1}{|c|}{\textbf{SUMMANDS}}&\multicolumn{1}{|c|}{\textbf{1}}
\\\hline
$ \omega_1'^{2}f_1\omega_1'\omega_2' \otimes
\omega_1'^{2}\omega_2'f_1 \otimes \omega_1'f_1\omega_2' \otimes
f_1\omega_2'\otimes f_2$ & $ax^{5}y$
\\\hline
$ \omega_1'^{3}\omega_2'f_1 \otimes \omega_1'^{2}f_1\omega_2'
\otimes \omega_1'f_1\omega_2' \otimes f_1\omega_2'\otimes f_2$ &
$ax^{6}y$
\\\hline
\end{tabular}
\end{center}


\begin{center}
 \large
\begin{tabular}{|l|p{2.5em}|}\hline
\multicolumn{1}{|c|}{\textbf{SUMMANDS}}&\multicolumn{1}{|c|}{\textbf{2}}
\\\hline
$ f_1\omega_1'^{3}\omega_2' \otimes f_1\omega_1'^{2}\omega_2'
\otimes f_1\omega_1'\omega_2' \otimes f_2\omega_1'\otimes f_1$ &
$a$
\\\hline
$ f_1\omega_1'^{3}\omega_2' \otimes \omega_1'f_1\omega_1'\omega_2'
\otimes f_1\omega_1'\omega_2' \otimes f_2\omega_1'\otimes f_1$ &
 $ax$
\\\hline
$ \omega_1'f_1\omega_1'^{2}\omega_2' \otimes
f_1\omega_1'^{2}\omega_2' \otimes f_1\omega_1'\omega_2' \otimes
f_2\omega_1'\otimes f_1$ &
 $ax$
\\\hline
$ \omega_1'^{2}f_1\omega_1'\omega_2' \otimes
f_1\omega_1'^{2}\omega_2' \otimes f_1\omega_1'\omega_2' \otimes
f_2\omega_1'\otimes f_1$ &
  $ax^2$
\\\hline

$ \omega_1'f_1\omega_1'^{2}\omega_2' \otimes
\omega_1'f_1\omega_1'\omega_2' \otimes f_1\omega_1'\omega_2'
\otimes f_2\omega_1'\otimes f_1$ &
  $ax^2$
\\\hline

$ \omega_1'^{2}f_1\omega_1'\omega_2' \otimes
\omega_1'f_1\omega_1'\omega_2' \otimes f_1\omega_1'\omega_2'
\otimes f_2\omega_1'\otimes f_1$ & $ax^3$
\\\hline
$ f_1\omega_1'^{3}\omega_2' \otimes f_1\omega_1'^{2}\omega_2'
\otimes \omega_1'\omega_2'f_1 \otimes \omega_1'f_2\otimes f_1$ &
$ax\bar{x}y$
\\
\hline
\end{tabular}
\end{center}

\begin{center}
 \large
\begin{tabular}{|l|p{2.5em}|}\hline
\multicolumn{1}{|c|}{\textbf{SUMMANDS}}&\multicolumn{1}{|c|}{\textbf{3}}
\\\hline
$ f_1\omega_1'^{3}\omega_2' \otimes f_1\omega_1'^2\omega_2'
\otimes \omega_1'f_2\omega_1' \otimes f_1\omega_1'\otimes f_1$ &
 $a\bar{x}$
\\\hline
$ f_1\omega_1'^{3}\omega_2' \otimes \omega_1'f_1\omega_1'\omega_2'
\otimes \omega_1'f_2\omega_1' \otimes f_1\omega_1'\otimes f_1$ &
$ax\bar{x}$
\\\hline
$ f_1\omega_1'^{3}\omega_2' \otimes \omega_1'^{2}\omega_2'f_1
\otimes \omega_1'^{2}f_2 \otimes f_1\omega_1'\otimes f_1$ &
  $ax^2\bar{x}^{2}y$
\\\hline

$ \omega_1'f_1\omega_1'^{2}\omega_2' \otimes
f_1\omega_1'^{2}\omega_2' \otimes \omega_1'f_2\omega_1' \otimes
f_1\omega_1'\otimes f_1$ &
 $ax\bar{x}$
\\\hline
$ \omega_1'^{2}f_1\omega_1'\omega_2' \otimes
f_1\omega_1'^{2}\omega_2' \otimes \omega_1'f_2\omega_1' \otimes
f_1\omega_1'\otimes f_1$ &
  $ax^2\bar{x}$
\\\hline

$ \omega_1'^{3}\omega_2'f_1 \otimes f_1\omega_1'^{2}\omega_2'
\otimes \omega_1'^{2}f_2 \otimes f_1\omega_1'\otimes f_1$ &
 $ax^3\bar{x}^{2}y$
\\\hline
$ \omega_1'f_1\omega_1'^{2}\omega_2' \otimes
\omega_1'f_1\omega_1'\omega_2' \otimes \omega_1'f_2\omega_1'
\otimes f_1\omega_1'\otimes f_1$ &
 $ax^2\bar{x}$
\\
\hline
\end{tabular}
\end{center}

\begin{center}
 \large
\begin{tabular}{|l|p{2.5em}|}\hline
\multicolumn{1}{|c|}{\textbf{SUMMANDS}}&\multicolumn{1}{|c|}{\textbf{3}}
\\\hline
$ f_1\omega_1'^{3}\omega_2' \otimes f_1\omega_1'^{2}\omega_2'
\otimes \omega_1'f_2\omega_1' \otimes \omega_1'f_1\otimes f_1$ &
 $ax\bar{x}$
\\\hline
$ f_1\omega_1'^{3}\omega_2' \otimes \omega_1'f_1\omega_1'\omega_2'
\otimes \omega_1'f_2\omega_1' \otimes \omega_1'f_1\otimes f_1$ &
  $ax^2\bar{x}$
\\\hline
$ f_1\omega_1'^{3}\omega_2' \otimes \omega_1'^{2}\omega_2'f_1
\otimes \omega_1'^{2}f_2 \otimes \omega_1'f_1\otimes f_1$ &
  $ax^3\bar{x}^{2}y$
\\\hline
$ \omega_1'f_1\omega_1'^{2}\omega_2' \otimes
f_1\omega_1'^{2}\omega_2' \otimes \omega_1'f_2\omega_1' \otimes
\omega_1'f_1\otimes f_1$ &
  $ax^2\bar{x}$
\\\hline
$ \omega_1'^{2}f_1\omega_1'\omega_2' \otimes
f_1\omega_1'^{2}\omega_2' \otimes \omega_1'f_2\omega_1' \otimes
\omega_1'f_1\otimes f_1$ &
  $ax^3\bar{x}$
\\\hline
$ \omega_1'^{3}\omega_2'f_1 \otimes f_1\omega_1'^{2}\omega_2'
\otimes \omega_1'^{2}f_2 \otimes \omega_1'f_1\otimes f_1$ &
  $ax^4\bar{x}^{2}y$
\\\hline
$ \omega_1'f_1\omega_1'^{2}\omega_2' \otimes
\omega_1'f_1\omega_1'\omega_2' \otimes \omega_1'f_2\omega_1'
\otimes \omega_1'f_1\otimes f_1$ &
  $ax^3\bar{x}$
\\
\hline
\end{tabular}
\end{center}


\begin{center}
 \large
\begin{tabular}{|l|p{2.5em}|}\hline
\multicolumn{1}{|c|}{\textbf{SUMMANDS}}&\multicolumn{1}{|c|}{\textbf{4}}
\\\hline
$ f_1\omega_1'^{3}\omega_2' \otimes \omega_1'^{2}f_2\omega_1'
\otimes f_1\omega_1'^2 \otimes f_1\omega_1'\otimes f_1$ &
  $a\bar{x}^2$
\\\hline
$ f_1\omega_1'^{3}\omega_2' \otimes \omega_1'^{2}f_2\omega_1'
\otimes \omega_1'f_1\omega_1' \otimes f_1\omega_1'\otimes f_1$ &
  $ax\bar{x}^2$
\\\hline

$ f_1\omega_1'^{3}\omega_2' \otimes \omega_1'^{2}f_2\omega_1'
\otimes \omega_1'^{2}f_1 \otimes f_1\omega_1'\otimes f_1$ &
 $ax^2\bar{x}^2$
\\\hline
$ \omega_1'f_1\omega_1'^{2}\omega_2' \otimes
\omega_1'^{2}f_2\omega_1' \otimes f_1\omega_1'^{2} \otimes
f_1\omega_1'\otimes f_1$ & $ax\bar{x}^{2}$
\\\hline
$ \omega_1'^{2}f_1\omega_1'\omega_2' \otimes
\omega_1'^{2}f_2\omega_1' \otimes f_1\omega_1'^{2} \otimes
f_1\omega_1'\otimes f_1$ & $ax^2\bar{x}^{2}$
\\\hline
$ \omega_1'^{3}\omega_2'f_1 \otimes \omega_1'^{3}f_2 \otimes
f_1\omega_1'^{2} \otimes f_1\omega_1'\otimes f_1$ &
 $ax^3\bar{x}^{3}y$
\\\hline
$ \omega_1'f_1\omega_1'^{2}\omega_2' \otimes
\omega_1'^{2}f_2\omega_1' \otimes \omega_1'f_1\omega_1' \otimes
f_1\omega_1'\otimes f_1$ &
  $ax^2\bar{x}^2$
\\
\hline
\end{tabular}
\end{center}


\begin{center}
 \large
\begin{tabular}{|l|p{2.5em}|}\hline
\multicolumn{1}{|c|}{\textbf{SUMMANDS}}&\multicolumn{1}{|c|}{\textbf{5}}
\\\hline
$ \omega_1'^{3}f_2\omega_1' \otimes f_1\omega_1'^{3} \otimes
f_1\omega_1'^2 \otimes f_1\omega_1'\otimes f_1$ &
 $a\bar{x}^3$
\\\hline
$ \omega_1'^{3}f_2\omega_1' \otimes f_1\omega_1'^{3} \otimes
\omega_1'f_1\omega_1' \otimes f_1\omega_1'\otimes f_1$ &
  $ax\bar{x}^3$
\\\hline
$ \omega_1'^{3}f_2\omega_1' \otimes f_1\omega_1'^{3} \otimes
\omega_1'^{2}f_1 \otimes f_1\omega_1'\otimes f_1$ &
$ax^2\bar{x}^3$
\\\hline
\end{tabular}
\end{center}


\begin{center}
 \large
\begin{tabular}{|l|p{2.5em}|}\hline
\multicolumn{1}{|c|}{\textbf{SUMMANDS}}&\multicolumn{1}{|c|}{\textbf{5}}
\\\hline
$ \omega_1'^{3}f_2\omega_1' \otimes \omega_1'f_1\omega_1'^{2}
\otimes f_1\omega_1'^{2} \otimes f_1\omega_1'\otimes f_1$ &
 $ax\bar{x}^{3}$
\\\hline

$ \omega_1'^{3}f_2\omega_1' \otimes \omega_1'^{2}f_1\omega_1'
\otimes f_1\omega_1'^{2} \otimes f_1\omega_1'\otimes f_1$ &
$ax^2\bar{x}^{3}$
\\\hline
$ \omega_1'^{3}f_2\omega_1' \otimes \omega_1'^{3}f_1 \otimes
f_1\omega_1'^{2} \otimes f_1\omega_1'\otimes f_1$ &
 $ax^3\bar{x}^{3}$
\\\hline
$ \omega_1'^{3}f_2\omega_1' \otimes \omega_1'f_1\omega_1'^2
\otimes \omega_1'f_1\omega_1' \otimes f_1\omega_1'\otimes f_1$ &
  $ax^2\bar{x}^3$
\\\hline

$ \omega_1'^{3}f_2\omega_1' \otimes \omega_1'f_1\omega_1'^2
\otimes \omega_1'^{2}f_1 \otimes f_1\omega_1'\otimes f_1$ &
  $ax^3\bar{x}^3$
\\\hline
$ \omega_1'^{3}f_2\omega_1' \otimes \omega_1'^{2}f_1\omega_1'
\otimes \omega_1'f_1\omega_1' \otimes f_1\omega_1'\otimes f_1$ &
  $ax^3\bar{x}^3$
\\\hline
$ \omega_1'^{3}f_2\omega_1' \otimes \omega_1'^{3}f_1 \otimes
\omega_1'f_1\omega_1' \otimes f_1\omega_1'\otimes f_1$ &
  $ax^4\bar{x}^3$
\\
\hline
\end{tabular}
\end{center}

The relevant terms and the results of calculations of
$\Delta^{(4)}(X)$ for $X = f_{1}f_{2}f_{1}^{3}$ in $(\ast)$ are as
follows:

\hspace*{\fill}TABULAR 5\hspace*{\fill}
\begin{center}
 \large
\begin{tabular}{|l|p{2.5em}|}\hline
\multicolumn{1}{|c|}{\textbf{SUMMANDS}}&\multicolumn{1}{|c|}{\textbf{1}}
\\\hline
$ f_1\omega_2'\omega_1'^3 \otimes \omega_2'f_1\omega_1'^2 \otimes
\omega_2'f_1\omega_1' \otimes \omega_2'f_1\otimes f_2$ &$ay^{3}$
\\\hline
$ f_1\omega_2'\omega_1'^3 \otimes \omega_2'f_1\omega_1'^2 \otimes
\omega_1'\omega_2'f_1 \otimes \omega_2'f_1\otimes f_2$ &$axy^{3}$
\\\hline
$ f_1\omega_2'\omega_1'^3 \otimes \omega_1'\omega_2'f_1\omega_1'
\otimes \omega_2'f_1\omega_1' \otimes \omega_2'f_1\otimes f_2$
&$axy^{3}$
\\\hline
$ f_1\omega_2'\omega_1'^3 \otimes \omega_1'^2\omega_2'f_1 \otimes
\omega_2'f_1\omega_1' \otimes \omega_2'f_1\otimes f_2$
&$ax^{2}y^{3}$
\\\hline
$ f_1\omega_2'\omega_1'^3 \otimes \omega_1'\omega_2'f_1\omega_1'
\otimes \omega_1'\omega_2'f_1 \otimes \omega_2'f_1\otimes f_2$
&$ax^{2}y^{3}$
\\\hline
$ f_1\omega_2'\omega_1'^3 \otimes \omega_1'^2\omega_2'f_1 \otimes
\omega_1'\omega_2'f_1 \otimes \omega_2'f_1\otimes f_2$
&$ax^{3}y^{3}$
\\\hline
$ \omega_1'\omega_2'f_1\omega_1'^2 \otimes f_1\omega_2'\omega_1'^2
\otimes \omega_2'f_1\omega_1' \otimes \omega_2'f_1\otimes f_2$
&$axy^3$
\\
\hline
\end{tabular}
\end{center}

\begin{center}
 \large
\begin{tabular}{|l|p{2.5em}|}\hline
\multicolumn{1}{|c|}{\textbf{SUMMANDS}}&\multicolumn{1}{|c|}{\textbf{1}}
\\\hline
$ \omega_1'^3\omega_2'f_1  \otimes \omega_1'\omega_2'f_1\omega_1'
\otimes f_1\omega_1'\omega_2' \otimes \omega_2'f_1\otimes f_2$
&$ax^{4}y^3$
\\\hline
$ \omega_1'^2\omega_2'f_1 \omega_1' \otimes
\omega_1'^2\omega_2'f_1 \otimes f_1\omega_1'\omega_2' \otimes
\omega_2'f_1\otimes f_2$ &$ax^{4}y^3$
\\\hline
$ \omega_1'^3\omega_2'f_1  \otimes \omega_1'^2\omega_2'f_1 \otimes
f_1\omega_1'\omega_2' \otimes \omega_2'f_1\otimes f_2$
&$ax^{5}y^3$
\\\hline$
\omega_1'\omega_2'f_1\omega_1'^2  \otimes
\omega_1'\omega_2'f_1\omega_1' \otimes \omega_1'\omega_2'f_1
\otimes f_1\omega_2'\otimes f_2$ &$ax^{3}y^3$
\\\hline

$ \omega_1'\omega_2'f_1\omega_1'^2  \otimes
\omega_1'^2\omega_2'f_1 \otimes \omega_1'\omega_2'f_1 \otimes
f_1\omega_2'\otimes f_2$ &$ax^{4}y^3$
\\\hline
$ \omega_1'^2\omega_2'f_1\omega_1'  \otimes
\omega_1'\omega_2'f_1\omega_1' \otimes \omega_1'\omega_2'f_1
\otimes f_1\omega_2'\otimes f_2$ &$ax^{4}y^3$
\\\hline
$ \omega_1'^3\omega_2'f_1  \otimes \omega_1'\omega_2'f_1\omega_1'
\otimes \omega_1'\omega_2'f_1 \otimes f_1\omega_2'\otimes f_2$
&$ax^{5}y^3$
\\\hline
$ \omega_1'^2\omega_2'f_1\omega_1'  \otimes
\omega_1'^2\omega_2'f_1 \otimes \omega_1'\omega_2'f_1 \otimes
f_1\omega_2'\otimes f_2$ &$ax^{5}y^3$
\\\hline
$ \omega_1'^3\omega_2'f_1  \otimes \omega_1'^2\omega_2'f_1 \otimes
\omega_1'\omega_2'f_1 \otimes f_1\omega_2'\otimes f_2$
&$ax^{6}y^3$
\\\hline

\end{tabular}
\end{center}


\begin{center}
\large
\begin{tabular}{|l|p{2.5em}|}\hline
\multicolumn{1}{|c|}{\textbf{SUMMANDS}}&\multicolumn{1}{|c|}{\textbf{2}}
\\\hline
$ f_1\omega_2'\omega_1'^3 \otimes \omega_2'f_1\omega_1'^2 \otimes
\omega_2'f_1\omega_1' \otimes f_2\omega_1'\otimes f_1$ & $ay^2$
\\\hline
$ f_1\omega_2'\omega_1'^3 \otimes \omega_2'f_1\omega_1'^2 \otimes
\omega_1'\omega_2'f_1 \otimes f_2\omega_1'\otimes f_1$  & $axy^2$
\\\hline
$ f_1\omega_2'\omega_1'^3 \otimes \omega_1'\omega_2'f_1\omega_1'
\otimes \omega_2'f_1\omega_1' \otimes f_2\omega_1'\otimes f_1$ &
$axy^2$
\\\hline
$ f_1\omega_2'\omega_1'^3 \otimes \omega_1'^2\omega_2'f_1 \otimes
\omega_2'f_1\omega_1' \otimes f_2\omega_1'\otimes f_1$ &
$ax^{2}y^2$
\\\hline
$ f_1\omega_2'\omega_1'^3 \otimes \omega_1'\omega_2'f_1\omega_1'
\otimes \omega_1'\omega_2'f_1 \otimes f_2\omega_1'\otimes f_1$ &
$ax^{2}y^2$
\\\hline
$ f_1\omega_2'\omega_1'^3 \otimes \omega_1'^2\omega_2'f_1 \otimes
\omega_1'\omega_2'f_1 \otimes f_2\omega_1'\otimes f_1$ &
$ax^{3}y^2$
\\\hline
$ \omega_1'\omega_2'f_1\omega_1'^2 \otimes f_1\omega_2'\omega_1'^2
\otimes \omega_2'f_1\omega_1' \otimes f_2\omega_1'\otimes f_1$
&$axy^2$
\\
\hline

\end{tabular}
\end{center}


\begin{center}
 \large
\begin{tabular}{|l|p{2.5em}|}\hline
\multicolumn{1}{|c|}{\textbf{SUMMANDS}}&\multicolumn{1}{|c|}{\textbf{3}}
\\\hline
$ f_1\omega_2'\omega_1'^3 \otimes \omega_2'f_1\omega_1'^2 \otimes
f_2\omega_1'^{2} \otimes f_1\omega_1'\otimes f_1$ & $ay$
\\\hline
$ f_1\omega_2'\omega_1'^3 \otimes \omega_2'f_1\omega_1'^2 \otimes
f_2\omega_1'^{2} \otimes \omega_1'f_1\otimes f_1$  & $axy$
\\\hline
$ f_1\omega_2'\omega_1'^3 \otimes \omega_1'\omega_2'f_1\omega_1'
\otimes f_2\omega_1'^{2} \otimes f_1\omega_1'\otimes f_1$ & $axy$
\\\hline

$ f_1\omega_2'\omega_1'^3 \otimes \omega_1'^2\omega_2'f_1 \otimes
f_2\omega_1'^{2} \otimes f_1\omega_1'\otimes f_1$  & $ax^{2}y$
\\\hline

$ f_1\omega_2'\omega_1'^3 \otimes \omega_1'\omega_2'f_1\omega_1'
\otimes f_2\omega_1'^{2} \otimes \omega_1'f_1\otimes f_1$ &
$ax^{2}y$
\\\hline
$ f_1\omega_2'\omega_1'^3 \otimes \omega_1'^2\omega_2'f_1 \otimes
f_2\omega_1'^{2} \otimes \omega_1'f_1\otimes f_1$& $ax^{3}y$
\\\hline
\end{tabular}
\end{center}

\begin{center}
 \large
\begin{tabular}{|l|p{2.5em}|}\hline
\multicolumn{1}{|c|}{\textbf{SUMMANDS}}&\multicolumn{1}{|c|}{\textbf{3}}
\\\hline
$ \omega_1'\omega_2'f_1\omega_1'^2 \otimes f_1\omega_2'\omega_1'^2
\otimes f_2\omega_1'^{2} \otimes f_1\omega_1'\otimes f_1$&$axy$
\\\hline
$ \omega_1'\omega_2'f_1\omega_1'^2 \otimes f_1\omega_2'\omega_1'^2
\otimes f_2\omega_1'^{2} \otimes \omega_1'f_1\otimes
f_1$&$ax^{2}y$
\\\hline
$ \omega_1'^2\omega_2'f_1\omega_1' \otimes f_1\omega_2'\omega_1'^2
\otimes f_2\omega_1'^{2} \otimes f_1\omega_1'\otimes
f_1$&$ax^{2}y$
\\\hline
$ \omega_1'^3\omega_2'f_1 \otimes f_1\omega_2'\omega_1'^2 \otimes
f_2\omega_1'^{2} \otimes f_1\omega_1'\otimes f_1$&$ax^{3}y$
\\\hline
$ \omega_1'^2\omega_2'f_1\omega_1' \otimes f_1\omega_2'\omega_1'^2
\otimes f_2\omega_1'^{2} \otimes \omega_1'f_1\otimes
f_1$&$ax^{3}y$
\\\hline
$ \omega_1'^3\omega_2'f_1 \otimes f_1\omega_2'\omega_1'^2 \otimes
f_2\omega_1'^{2} \otimes \omega_1'f_1\otimes f_1$ &$ax^{4}y$
\\\hline
$ \omega_1'\omega_2'f_1\omega_1'^2  \otimes
\omega_1'\omega_2'f_1\omega_1' \otimes \omega_1'f_2\omega_1'
\otimes f_1\omega_1'\otimes f_1$ &$ax^{2}\bar{x}y^2$
\\
\hline
\end{tabular}
\end{center}


\begin{center}
\large
\begin{tabular}{|l|p{2.5em}|}\hline
\multicolumn{1}{|c|}{\textbf{SUMMANDS}}&\multicolumn{1}{|c|}{\textbf{4}}
\\\hline
$ f_1\omega_2'\omega_1'^3 \otimes f_2\omega_1'^3 \otimes
f_1\omega_1'^2 \otimes f_1\omega_1'\otimes f_1$  & $a$
\\\hline
$ f_1\omega_2'\omega_1'^3 \otimes f_2\omega_1'^3 \otimes
f_1\omega_1'^2 \otimes \omega_1'f_1\otimes f_1$ & $ax$
\\\hline
$ f_1\omega_2'\omega_1'^3 \otimes f_2\omega_1'^3 \otimes
\omega_1'f_1\omega_1'\otimes f_1\omega_1'\otimes f_1$ & $ax$
\\\hline
$ f_1\omega_2'\omega_1'^3 \otimes f_2\omega_1'^3 \otimes
\omega_1'^{2}f_1\otimes f_1\omega_1'\otimes f_1$ & $ax^2$
\\\hline

$ f_1\omega_2'\omega_1'^3 \otimes f_2\omega_1'^3 \otimes
\omega_1'f_1\omega_1' \otimes \omega_1'f_1\otimes f_1$ & $ax^2$
\\\hline
$ f_1\omega_2'\omega_1'^3 \otimes f_2\omega_1'^3 \otimes
\omega_1'^{2}f_1 \otimes \omega_1'f_1\otimes f_1$  & $ax^3$
\\\hline
$ \omega_1'\omega_2'f_1\omega_1'^2 \otimes
\omega_1'f_2\omega_1'^{2} \otimes f_1\omega_1'^2 \otimes
f_1\omega_1'\otimes f_1$  & $ax\bar{x}y$
\\
\hline
\end{tabular}
\end{center}

\begin{center}
\large
\begin{tabular}{|l|p{2.5em}|}\hline
\multicolumn{1}{|c|}{\textbf{SUMMANDS}}&\multicolumn{1}{|c|}{\textbf{4}}
\\\hline
$ \omega_1'^2\omega_2'f_1\omega_1'  \otimes
\omega_1'f_2\omega_1'^{2} \otimes \omega_1'^{2}f_1 \otimes
\omega_1'f_1\otimes f_1$ &$ax^{5}\bar{x}y$
\\\hline
$ \omega_1'^3\omega_2'f_1  \otimes \omega_1'f_2\omega_1'^{2}
\otimes \omega_1'^{2}f_1 \otimes \omega_1'f_1\otimes
f_1$&$ax^{6}\bar{x}y$
\\\hline
\end{tabular}
\end{center}


\begin{center}
\large
\begin{tabular}{|l|p{2.5em}|}\hline
\multicolumn{1}{|c|}{\textbf{SUMMANDS}}&\multicolumn{1}{|c|}{\textbf{5}}
\\\hline
$ \omega_1'f_2\omega_1'^3 \otimes f_1\omega_1'^3 \otimes
f_1\omega_1'^2 \otimes f_1\omega_1'\otimes f_1$& $a\bar{x}$
\\\hline
$ \omega_1'f_2\omega_1'^3 \otimes f_1\omega_1'^3 \otimes
f_1\omega_1'^2 \otimes \omega_1'f_1\otimes f_1$ & $ax\bar{x}$
\\\hline
$ \omega_1'f_2\omega_1'^3 \otimes f_1\omega_1'^3 \otimes
\omega_1'f_1\omega_1' \otimes f_1\omega_1'\otimes f_1$ &
$ax\bar{x}$
\\\hline
$ \omega_1'f_2\omega_1'^3 \otimes f_1\omega_1'^3 \otimes
\omega_1'^{2}f_1 \otimes f_1\omega_1'\otimes f_1$& $ax^2\bar{x}$
\\\hline
$ \omega_1'f_2\omega_1'^3 \otimes f_1\omega_1'^3 \otimes
\omega_1'f_1\omega_1' \otimes \omega_1'f_1\otimes f_1$
&$ax^2\bar{x}$
\\\hline
$ \omega_1'f_2\omega_1'^3 \otimes f_1\omega_1'^3 \otimes
\omega_1'^{2}f_1 \otimes \omega_1'f_1\otimes f_1$& $ax^3\bar{x}$
\\\hline
$ \omega_1'f_2\omega_1'^3 \otimes \omega_1'f_1\omega_1'^2 \otimes
f_1\omega_1'^2 \otimes f_1\omega_1'\otimes f_1$&$ax\bar{x}$
\\
\hline
\end{tabular}
\end{center}

\bigskip
\bibliographystyle{amsalpha}

\begin{thebibliography}{A}
\medskip



\bibitem [BGH1]{BGH} N. Bergeron, Y. Gao and N. Hu, \textit{Drinfel'd doubles and
Lusztig's symmetries of two-parameter quantum groups},
math.RT/0505614, Journal of Algebra (to appear), in 2006.

\bibitem [BGH2]{BGH} N. Bergeron, Y. Gao and N. Hu, \textit{Representations of
two-parameter quantum orthogonal groups and symplectic groups},
math.QA/0510124, Contemp. Math. of AMS (to appear), in 2006.

\bibitem [BKL]{BKL} G. Benkart, S.~J. Kang and K.~H. Lee, \textit{On the center of two-parameter quantum groups},
Preprint 2004.

\bibitem [BW1]{B} G.
Benkart and S. Witherspoon, \textit{Two-parameter quantum groups
and Drinfel'd doubles}, Algebr. Represent. Theory, \textbf{7}
(2004), 261--286.

\bibitem [BW2]{B} G. Benkart and S. Witherspoon, \textit{Representations of two-parameter quantum
groups and Schur-Weyl duality}, Hopf algebras, pp. 65--92, Lecture
Notes in Pure and Appl. Math., 237, Dekker, New York, 2004.

\bibitem [BW3]{B} G. Benkart and S. Witherspoon, \textit{Restricted two-parameter quantum groups},
Fields Institute Communications, ``Representations of Finite
Dimensional Algebras and Related Topics in Lie Theory and
Geometry", vol. 40, Amer. Math. Soc., Providence, RI, 2004, pp.
293--318.

\bibitem [Dr]{Dr} V.G. Drinfel'd, \textit{Quantum groups}, in
``Proceedings ICM'', Berkeley, Amer. Math. Soc. 1987, pp.
798--820.

\bibitem [KS]{KS} A. Klimyk and K. Schm\"udgen, \textit{Quantum Groups and
Their Representations}, Springer-Verlag Berlin Heidelberg New
York, 1997.

\bibitem [L1]{L} G. Lusztig, \textit{Introduction to Quantum Groups},
Birkh\"auser Boston, 1993.

\bibitem [L2]{L} G. Lusztig, \textit{Quantum groups at roots of $1$}, Geom. Dedicata \textbf{35}
(1990), 89--114.

\bibitem [RTF]{RTF} N.Yu. Reshetikhin, L.A. Takhtajan and L.D. Faddeev,
\textit{Quantization of Lie groups and Lie algebras}, Algebra and
Anal. \textbf{1}, 178--206 (1989) (Leningrad Math. J. 1
     [Engl. transl. 193--225 (1990)]).
\end{thebibliography}

\end{document}